\documentclass[11pt]{amsart}
\usepackage{amssymb}
\newfont{\ams}{msbm10 at 12pt}
\newfont{\amsi}{msbm8}
\newcommand{\nc}{\newcommand}
\nc{\C}{{\mathcal C}}
\nc{\CB}{{\mathcal B}}
\nc{\supp}{\operatorname{supp}}
\nc{\M}{{\mathcal M}}
\nc{\LL}{{\mathbf L}}
\renewcommand{\P}{{\mathbf  P}}
\renewcommand{\L}{{\mathcal L}}
\nc{\D}{{\mathbb D}}
\renewcommand{\O}{{\mathcal O}}
\nc{\St}{\operatorname{St}^{\bullet}}
\nc{\A}{{\frak A}}

\nc{\1}{{\mathbf 1}}
\nc{\CC}{{\mathbb C}}
\nc{\R}{{\mathbb R}}
\renewcommand{\k}{{\mathbb  k}}
\nc{\Q}{{\mathbb Q}}
\nc{\U}{{\mathbf U}}
\nc{\B}{{\mathbf B}}
\nc{\N}{{\mathbb N}}
\nc{\Z}{{\mathbb Z}}
\nc{\Rhom}{\operatorname{RHom}\bul}
\nc{\Ad}{\operatorname{Ad}}
\nc{\Res}{\operatorname{Res}}
\nc{\gr}{\operatorname{gr}}
\nc{\tr}{\operatorname{tr}}
\nc{\End}{\operatorname{End}}

\nc{\g}{{\frak g}}
\nc{\hatg}{\hat{\frak g}}

\renewcommand{\b}{{\frak b}}
\nc{\sem}{{$\S_{\g}^{\g_{>0}}}}
\nc{\gl}{{\frak g\frak l}}
\nc{\n}{{\frak n}}
\nc{\si}{{\frac\infty 2}}
\nc{\p}{{\frak p}}
\nc{\h}{{\frak h}}
\renewcommand{\u}{{\frak u}}
\nc{\Ind}{\operatorname{Ind}}
\nc{\ch}{\operatorname{ch}}
\nc{\Coind}{\operatorname{Coind}}
\nc{\opp}{{\operatorname{opp}}}
\nc{\Ker}{\operatorname{Ker}}
\nc{\im}{\operatorname{Im}}
\nc{\Coker}{\operatorname{Coker}}
\nc{\dirlim}{\underset{\rightarrow}{\operatorname{lim}}}
\nc{\invlim}{\underset{\leftarrow}{\operatorname{lim}}}
\nc{\Sem}{{{\mathbf S}_{\g}^{\g_{>0}}}}
\nc{\CN}{{\mathcal N}}
\nc{\Ext}{\operatorname{Ext}^{\bullet}}
\nc{\ext}{\operatorname{Ext}}
\nc{\tilW}{\til{W}}

\nc{\lth}{\ell t}
\nc{\BB}{\mathcal{B}}
\nc{\Tor}{\operatorname{Tor}_{\bullet}}
\nc{\tor}{\operatorname{Tor}}
\nc{\Tors}{\operatorname{Tor}_{\frac \infty 2+\bullet}}
\nc{\Exts}{\operatorname{Ext}^{\frac \infty 2+\bullet}}
\nc{\Hom}{\operatorname{Hom}^{\bullet}}
\nc{\ad}{\operatorname{ad}}
\renewcommand{\hom}{\operatorname{Hom}}

\renewcommand{\mod}{\operatorname{-mod}}
\nc{\Mod}{\operatorname{Mod}}

\nc{\Barb}{\operatorname{Bar}^{\bullet}}

\nc{\upX}{X^{\uparrow}}
\nc{\upcD}{{\mathcal D}^{\uparrow}}
\nc{\upD}{D^{\uparrow}}
\nc{\dX}{X^{\downarrow}}
\nc{\dcD}{{\mathcal D}^{\downarrow}}
\nc{\dD}{D^{\downarrow}}
\nc{\upC}{{\mathcal C}^{\uparrow}}
\nc{\dC}{{\mathcal C}^{\downarrow}}
\nc{\underA}{\underline{A}}
\nc{\underC}{\underline{\CC}}
\nc{\underB}{\underline{B}}
\nc{\underk}{\underline{\k}}
\nc{\Db}{D^{\bullet}}
\nc{\ten}{{\otimes}}
\nc{\tenl}{\overset{\operatorname{L}}\ten}
\nc{\map}{\longrightarrow}

\nc{\bs}{\bigskip\\}
\nc{\ms}{\smallskip\\}

\nc{\tilbar}{\widetilde{\operatorname{Bar}}}
\nc{\tilBarb}{\widetilde{\operatorname{Bar}}^{\bullet}}
\nc{\overr}{\overline{R}}
\nc{\overI}{\overline{I}}
\nc{\overX}{\overline{X}}
\nc{\overh}{\overline{h}}
\nc{\overY}{\overline{Y}}
\nc{\overW}{\overline{W}}
\nc{\linbar}{\overline{\operatorname{Bar}}}

\nc{\til}{\widetilde}
\nc{\oppA}{A^{\sharp}}
\nc{\Lemma}{{\bf Lemma:\ }}
\nc{\Theorem}{{\bf Theorem:}\ }
\nc{\Cor}{{\bf Corollary:}\ }
\nc{\Def}{{\bf Definition:}\ }
\nc{\Prop}{{\bf Proposition:}\ }
\nc{\Con}{{\bf Conjecture:}\ }
\nc{\Rem}{{\bf Remark:}\ }
\nc{\dok}{{\bf Proof.}\ }

\nc{\SInd}{\operatorname{S-Ind}}
\nc{\SCoind}{\operatorname{S-Coind}}
\nc{\bul}{^{\bullet}}
\nc{\stand}{C^{\frac\infty2+\bullet}}
\nc{\ssn}{\subsection{}}
\nc{\sssn}{\subsubsection{}}

\nc{\hgt}{\operatorname{ht}}
\nc{\sqbinom}{\fracwithdelims[][0pt]}

\address{Independent University of Moscow, Pervomajskaya st. 16-18,
Moscow 105037, Russia}
\email{hippie@@mccme.ru}
\author{Sergey Arkhipov}
\title{A proof of Feigin's conjecture}
\date{}

\begin{document}
\maketitle
\section{Introduction}
     The present paper is devoted to further development of
semiinfinite cohomology of small quantum groups. The topic
appeared first in \cite{Ar1}, where the very definition of
semiinfinite cohomology $\Exts_A(\underC,\cdot)$ was given. The
setup for the definition of semiinfinite cohomology of an
algebra $A$ includes two subalgebras $B,N\subset A$ and the
{\em triangular decomposition} of $A$, i. e. the vector space
isomorphism $B\ten N\til{\map}A$ provided by the multiplication
in $A$. Fix root data $(Y,X,\ldots)$ of the finite type
$(I,\cdot)$ and a positive integer number $\ell$. The {\em
small quantum group} $\u_\ell$ with the standard triangular
decomposition turns out to be a very interesting object for the
investigation of semiinfinite cohomology. The explanation for
this lies in the following fact proved by Ginzburg and Kumar in
\cite{GK}. Consider the set of nilpotent elements $\CN$ in the
simple Lie algebra $\g$ corresponding to $(Y,X,\ldots)$. \vskip
1mm \noindent {\bf Theorem:}
$\Ext_{\u_\ell}(\underC,\underC)={\mathcal F}un(\CN)$ as an
associative algebra. The grading on the right hand side is
provided by the action of the group $\CC^*$ on the affine
variety $\CN$.\qed

     On the other hand it is proved in \cite{Ar2} that the
algebra $\Ext_A(\underC,\underC)$ acts naturally on the
semiinfinite cohomology of $A$. Thus in particular
$\Exts_{\u_\ell}(\underC,\cdot)$ becomes a quasicoherent sheaf
on $\CN$. It is natural to look for the answer for semiinfinite
cohoomlogy of $\u_\ell$ in terms of geometry of $\CN$. B.
Feigin has proposed the following conjecture. Consider the
standard positive nilpotent subalgebra
$\n^+\subset\CN\subset\g$.
\vskip 3mm
\noindent
{\bf Conjecture A:}
The quasicoherent sheaf on $\CN$ provided by
$\Exts_{\u_\ell}(\underC,\underC)$ is equal to the sheaf of algebraic
distributions on $\CN$ with support on $\n^+\subset\CN$.\qed

Moreover note that the simply connected Lie group $G$ with the Lie
algebra equal to $\g$ acts naturally on $\CN$. This action provides a
structure of a $\n^+$-integrable $\g$-module on the described
distributions' space.  On the other hand it was shown in \cite{Ar1}
that there exists a natural $U(\g)$-module structure on
$\Exts_{\u_\ell}(\underC,\underC)$.  The $\g$-module version of the
Feigin conjecture states that the described $\g$-modules are
isomorphic.  In \cite{Ar1} and \cite{Ar5} the conjecture was proved
on the level of characters of the $\g$-modules.  In the present paper
we give a full proof of the Feigin conjecture.

\ssn
Let us describe briefly the structure of the paper. In the second
section we recall nessesary facts concerning the usual cohomology of
small quantum groups. The main facts here include the mentioned
Ginzburg-Kumar theorem, the Kostant theorem describing the algebra
structure on ${\mathcal F}un(\CN)$ (see \ref{kost}) and a homological
description of a certain degeneration $\til{\u}_\ell$ of the algebra
$\u_\ell$ (see \ref{degen}). It turns out that both the usual and the
semiinfinite cohoomlogy of the algebra $\til{\u}_\ell$ with trivial
coefficients can be described as associated graded objects for the
corresponding functors over $\u_\ell$ with respect to  certain
filtrations.

     In the third section we prove the ${\mathcal
Fun}(\CN)$-module version of the Feigin conjecture. The main
arguments here are as follows. First we describe the
$\Ext_{\til{\u}_\ell}(\underC,\underC)$-module
$\Exts_{\til{\u}_\ell}(\underC,\underC)$. Next, using this
description, we construct a morphism of ${\mathcal
F}un(\CN)$-modules
$$
\Phi:\
H^0(\CN\setminus D_1\cup\ldots\cup D_{\sharp(R^+)},\O_{\CN})
\map\Exts_{\u_\ell}(\underC,\underC),
$$
     where $D_1,\ldots,D_{\sharp(R^+)}$ denote the coordinate
divisors such that $\cap_i D_i=\n^+$. Note that the first
${\mathcal F}un(\CN)$-module contains the distributions module
$H^{\sharp(R^+)}_{\n^+}(\CN,\O_{\CN})$ as a certain quotient
module.

     Finally, using the connection between
$\Exts_{\u_\ell}(\underC,\underC)$
     and $\Exts_{\til{\u}_\ell}(\underC,\underC)$ we prove that
the map $\Phi$ provides a ${\mathcal F}un(\CN)$-module
isomorphism $$H^{\sharp(R^+)}_{\n^+}(\CN,\O_{\CN})\til{\map}
\Exts_{\u_\ell}(\underC,\underC).$$ In the fourth section we
prove the $\g$-module version of the Feigin conjecture. The
main steps of the proof are as follows. First we construct
explicitly the $\n^+$-module isomorphism. Next we recall from
\cite{Ar1} that there exists a nondegenerate $\g$-equivariant
contragradient pairing on $\Exts_{\u_\ell}(\underC,\underC)$.
We construct its geometric analogue on
$H^{\sharp(R^+)}_{\n^+}(\CN,\O_{\CN})$. On the other hand it is
easy to verify that $H^{\sharp(R^+)}_{\n^+}(\CN,\O_{\CN})$ is
free over the algebra $U(\n^-)$. It follows that both
$H^{\sharp(R^+)}_{\n^+}(\CN,\O_{\CN})$ and
$H^{\sharp(R^+)}_{\n^+}(\CN,\O_{\CN})$ are co-free over
$U(\n^+)$. Using the constructed contragradient pairung on the
$\g$-module of semiinfinite cohomology we see that the latter
module is $U(\n^-)$-free. It follows that both
$H^{\sharp(R^+)}_{\n^+}(\CN,\O_{\CN})$ and
$H^{\sharp(R^+)}_{\n^+}(\CN,\O_{\CN})$ are {\em tilting}
$\g$-modules. Finally a beautiful result of Andersen states
that a tilting $\g$-module is completely determined up to an
isomorphism by its character (see \cite{A}).

     In the fifth section we p resent a generalization of the
Feigin conjecture as follows. Consider the {\em contragradient
Weyl module} $\D W(\ell\lambda)$ over $\u_\ell$ with the
highest weight $\ell\lambda$. On the other hand let $p:\
T^*(G/B)\map G/B)$ denote the cotangent bundle to the flag
variety of $G$. The moment map for the symplectic $G$-action on
$T^*(G/B)$ provides the {\em Springer-Grothendieck resolution}
of singularities of the nilpotent cone $\mu:\ T^*(G/B)\map\CN$.
Consider the linear bundle $\L(\lambda)$ on $G/B$. \vskip 3mm
\noindent {\bf Conjecture B:} The ${\mathcal
F}un(\CN)=\Ext_{\u_\ell}(\underC,\underC)$-module
$\Exts_{\u_\ell}(\underC,\D W(\ell\lambda))$ is isomorphic to
$H^{\sharp(R^+)}_{\n^+}(\CN,\mu_*p^*\L(\lambda))$.\qed

     We give a sketch of the proof of the Conjecture B. The
main tool here is a certain specialization of the quantum BGG
resolution into the root of unity. We call it {\em the
contragradient quasi-BGG complex}. This complex has $\D
W(\ell\lambda)$ as a zero cohomology module. We conjecture that
the contragradient quasi-BGG complex is in fact quasiisomorphic
to $\D W(\ell\lambda)$. Still even without the last assumption
we manage to prove Conjecture B.

\subsection{Acknowledgements.}
     The material of the paper is based upon work supported by
the U.S. Civilian Research and Development Foundation under
Award No. RM1-265 and by the grant INTAS 94-94720.

     Part of the work over results of the paper was done during
the author's visit to Ecole Polytechnique in Paris in December
1996. The author is happy to thank the inviting orgainization
for hospitality and stimulating research conditions.

     The author would like to thank B. Feigin, M. Finkelberg,
V. Ginzburg, W. Soergel and V. Ostrik for helpful discussions.

\section{Small quantum groups.}
\subsection{Root data.} Fix a  {\em Cartan datum}
$(I,\cdot)$ of the finite type   and a {\em simply connected root
datum} $(Y,X,\ldots)$ of the type $(I,\cdot)$. Thus we have
$Y=\Z[I]$, $X=\hom(Y,\Z)$, and the pairing $\langle\ ,\ \rangle:\
Y\times X\map\Z$ coincides
 with the natural one (see \cite{L3},  I 1.1, I
2.2). In particular the data contain canonical embeddings
$I\hookrightarrow Y, i\mapsto i$ and $I\hookrightarrow X, i\mapsto i':
\langle i',j\rangle:= 2i\cdot j/i\cdot i$. The latter map is naturally
extended to an embedding $Y\subset X$. Denote by $\hgt$ the linear
function on  $X$ defined on elements $i', i\in I$, by $\hgt(i')=1$ and
extended to the whole $X$ by linearity. The root system (resp. the
positive root system) corresponding to the data $(Y,X,\ldots)$ is
denoted by $R$ (resp. by $R^+$), below $W$ denotes the Weyl group
of $R$.

\subsection{Quantum groups at roots of 1.}
     Given the root datum $(Y,X,\ldots)$ Drinfeld and Jimbo
constructed an associative algebra $\U$ over the field $\Q(v)$
of rational functions in $v$ with generators $E_i, F_i,
K_i^{\pm1}, i\in I$, and relations being the quantum analogues
of the classical Serre relations in the universal enveloping
algebra of the corresponding simple Lie algebra $\g$. The
explicit form of the relations can be found e.~g. in \cite{L1},
I 3.1. We call the algebras $\U$ the {\em quantum groups}.

     Lusztig (see \cite{L1}, V 31.1) defined a
$\Q[v,v^{-1}]$-subalgebra $\U_{\Q}$ in $\U$ being the quantum
analogue of the integral form for the universal enveloping
algebra of $\g$ due to Kostant. In particular the elements
$E_i, F_i, K_i^{\pm1},i\in I$, belong to $\U_{\Q}$. Let $\ell$
be an odd number satisfying the conditions from \cite{GK}. Fix
a primitive $\ell$-th root of unity $\zeta$. Define a
$\CC$-algebra $\til{\U}_\ell:=\U_{\Q}\ten_{\Q[v,v^{-1}]}\CC$,
where $v$ acts on $\CC$ by multiplication by $\zeta$. It is
known that the elements $K_i^{\ell}, i\in I$, are central in
$\til{\U}_\ell$. Set $\U_\ell:=\til{\U}_\ell/(K_i^\ell-1, i\in
I)$. The algebra $\U_\ell$ is generated by the elements $E_i$,
$E_i^{(\ell)}$, $F_i$, $F_i^{(\ell)}$, $K_i^{\pm1}$, $i\in I$.
Here $E_i^{(\ell)}$ (resp. $F_i^{(\ell)}$) denotes the {\em
$\ell$-th quantum divided power} of the element $E_i$ (resp.
$F_i$) specialized at the ro ot of unity $\zeta$.

Following Lusztig we
define the {\em small quantum group} $\u_\ell$ at the root of unity
$\zeta$ as the subalgebra in  $\U_\ell$ generated by all $E_i, F_i,
K_i^{\pm1}, i\in I$. Denote the subalgebra in  $\u_\ell$ generated by
$E_i, i\in I$ (resp. $F_i, i\in I$, resp. $K_i, i\in I$), by
$\u_\ell^+$ (resp.  $\u_\ell^-$, resp. $\u_\ell^0$). Note that the
algebra $\u_\ell$ is graded naturally by the abelian group  $X$. Using
the function $\hgt$ we obtain a $\Z$-grading on  $\u_\ell$ from this
$X$-grading. In particular the subalgebra $\u_\ell^+$ (resp.
$\u_\ell^-$) is graded by $\Z_{\ge0}$ (resp. by $\Z_{\le0}$).

     Below we present several well known facts about the
algebra $\u_\ell$ to be used later. Recall that an augmented
subalgebra $B\subset A$ with the augmentation ideal
$\overline{B}\subset B$ is called {\em normal} if
$A\overline{B}=\overline{B}A$. If so, the space
$A/A\overline{B}$ becomes an algebra. It is denoted by $A//B$.
Fix an augmentation on $\u_\ell$ as follows: $E_i\mapsto0,
F_i\mapsto 0, K_i\mapsto 1$ for every $i\in I$. Set
$\underC:=\u_\ell/\overline{\u}_\ell$.

 \sssn \Lemma (see \cite{AJS} 1.3, \cite{L2} Theorem 8.10)
\label{mainu} \begin{itemize} \item[(i)] The multiplication in
$\u_\ell$ provides a vector space isomorphism
 $\u_\ell=\u_\ell^-\ten\u_\ell^0\ten\u_\ell^+$;
$\dim\u_\ell^{\pm}=\ell^{\sharp(R^+)}$; the subalgebra $\u^0_\ell$
is isomorphic to the group algebra of the group $(\Z/\ell\Z)^{\sharp(I)}$.
\item[(ii)] The subalgebra $\u_\ell\subset \U_\ell$ is normal and we have
$\U_\ell//\u_\ell=U(\g)$.\qed
\end{itemize}
The subalgebra $\u_\ell^-\ten\u_\ell^0$ (resp.
$\u_\ell^0\ten\u_\ell^+$) in  $\u_\ell$ is denoted by $\b_\ell^-$
(resp. by $\b_\ell^+$). Recall that a finite dimensional algebra $A$
is called {\em Frobenius} if the left  $A$-modules $A$ and
$A^*:=\hom_{\CC}(A,\CC)$ are isomorphic.

\sssn
\Lemma  (see \cite{Ar1}, Lemma 2.4.5)                                          \label{one}
The algebras  $\u_\ell^+$  and  $\u_\ell^-$ are Frobenius.\qed

Consi
der the filtration on the algebra  $\u_\ell$ as follows. Let the
filtration component $F^{\le d}(\u_\ell)$ be linearly generated by
$X$-homogeneous monomials $u=u^-\ten u^0\ten u^+$ such that
$|\hgt(\deg u^-)|+|\hgt(\deg u^+)|\le d$. By definition set
$\til{\u}_\ell:=\gr^{F}\u_\ell$.  Evidently we have
$\gr^F\b_\ell^+=\b_\ell^+$, $\gr^F\b_\ell^-=\b_\ell^-$,
$\gr^F\u_\ell^0=\u_\ell^0$.

\sssn
\Lemma                            \label{commute}
Elements of the subalgebra $\u_\ell^-\subset\til{\u}_\ell$ commute with
elements of the subalgerba $\u_\ell^+\subset\til{\u}_\ell$.\qed

Both the $X$-  and the  $\Z$-grading as well as the augmentation on
$\til{\u}_\ell$ are induced by the ones on  $\u_\ell$.  Denote the
category of $X$-graded finite dimensional left $\u_\ell$-modules
(resp. $\til{\u}_\ell$-modules) $M=\underset{\lambda\in
X}{\bigoplus}M_\lambda$ such that $K_i$ acts on  $M_\lambda$ by
multiplication by the scalar $\zeta^{\langle i,\lambda\rangle}$ and
$E_i:\ M_\lambda\map M_{\lambda+i'},\ F_i:\ M_\lambda\map
M_{\lambda-i'}$ for all $i\in I$,  with morphisms preserving
$X$-gradings, by $\u_\ell\mod$ (resp. by $\til{\u}_\ell\mod$).  For
$M,N\in \u_\ell\mod$ and $\lambda\in\ell\cdot X$ we define the
shifted module
$M\langle\lambda\rangle\in\u_\ell\mod:\
M\langle\lambda\rangle_\mu:=M_{\lambda+\mu}$ and set
$\hom_{\u_\ell}(M,N):=\underset{\lambda\in\ell\cdot
X}{\bigoplus}\hom_{\u_\ell\mod}(M\langle \lambda\rangle,N)$.  The
functor $\hom_{\til{\u}_\ell}$ is defined in a similar way. Evidently
the spaces $\hom_{\u_\ell}(\cdot,\cdot)$  and
$\hom_{\til{\u}_\ell}(\cdot,\cdot)$ posess  natural $\ell\cdot
X$-gradings.

     \subsection{Cohomology of small quantum groups.} Consider
the $\ell\cdot X\times\Z$-graded algebra
$\Ext_{\u_\ell}(\underC,\underC)$. Note that by Shapiro lemma
and Lemma~\ref{mainu} (ii) the Lie algebra $\g$ acts naturally
on the $\ext$ algebra and the multiplication in the algebra
satisfies Lebnitz rule with respect to the $\g$-action. In
\cite{GK} Ginzburg and Kumar obtained a nice description of the
mu ltiplication structure as well as the $\g$-module structure
on $\Ext_{\u_\ell}(\underC,\underC)$ as follows.

     \subsubsection{Functions on the nilpotent cone.} Let $G$
be the simply connected Lie group with the Lie algebra $\g$.
Then $G$ acts on $\g$ by adjunction. The action preserves the
set of nilpotent elements $\CN\subset\g$ called the {\em
nilpotent cone} of $\g$. The action is algebraic, thus it
provides a morphism of $\g$ into the Lie algebra of algebraic
vector fields on the nilpotent cone ${\mathcal V}ect(\CN)$. The
latter algebra acts on the algebraic functions
$H^0(\CN,\O_{\CN})$. The action is $G$-integrable. Note also
that the natural action of the group $\CC^*$ provides a grading
on $H^0(\CN,\O_{\CN})$ preserved by the $G$-action.

\sssn
\Theorem
     (see \cite{GK}) The algebra and $\g$-module structures on
$H^0(\CN,\O_{\CN})$ and on $\Ext_{\u_\ell}(\underC,\underC)$
coincide. The homological grading on the latter algebra
corresponds to the grading on the former one provided by the
$\CC^*$-action. The $X$-grading on $H^0(\CN,\O_{\CN})$ provided
by the weight decomposition with respect to the action of the
Cartan subalgebra in $\g$ corresponds to the natural $\ell\cdot
X$-grading on the space $\Exts_{\u_\ell}(\underC,\underC)$.
\qed

We will need a more detailed description of the algebra $H^0(\CN,\O_{\CN})$.
First the folowing result of Kostant shows the size of the algebra.

\sssn
\Prop (see \cite{H})
 $
 \ch\left(H^0(\CN,{\mathcal O}_{\CN}),t\right)=
{\frac{\underset{w\in W}{\sum}t^{l(w)}}
{\underset{\alpha\in R^+}{\prod}(1-e^{-\alpha}t)(1-e^{\alpha}t)}}$.

     Here the indeterminate $t$ stands for the homogeneous
$\Z$-grading, for a weight $\alpha=\sum a_ii'$ the symbol
$e^\alpha$ denotes the monomial $\prod_i(e^{i'})^{a_i}$, and
$l(w)$ denotes the length of the element of the Weyl group.
\qed

     In particular the space $$
\left(H^0(\CN,\O_{\CN})\right)_2= \underset{\alpha\in
R^+}{\bigoplus}\CC e_\alpha^*\oplus\underset{w\in
W,\lth(w)=1}{\bigoplus} \CC \xi_w \oplus\underset{\alpha\in
R^+}{\bigoplus} \CC f_\alpha^* =\n^{+*}\oplus \underset{w\in
W,l(w)=1}{\bigoplus}\CC\xi_w\oplus\n^{-*}, $$ where the
$X$-grading of the element $e_\alpha^*$ (resp. $f^*_\alpha$,
resp. $\xi_w$) is equal to $-\alpha$ (resp. to $\alpha$, resp.
to $0$). Consider the subalgebra $H_-$ (resp. $H_+$) in
$H^0(\CN,\O_{\CN})$ grnerated by the space $\n^{+*}$ (resp.
$\n^{-*}$). The detailed description of the algebra structure
on $H^0(\CN,\O_{\CN})$ is given by the following statement.

     \sssn \label{kost} \Prop (see \cite{Ko}, Theorem 1.5) The
algebra $H_-$ (resp. $H_+$) is equal to the free commutative
algebra $S\bul(\n^{+*})$ (resp. $S\bul(\n^{-*})$). \qed

     Note that the inclusion $S\bul(\n^{-*})\hookrightarrow
H^0(\CN,\O_{\CN})$ (resp. $S\bul(\n^{+*})\hookrightarrow
H^0(\CN,\O_{\CN})$) corresponds to the coordinate projection
$\CN\hookrightarrow\g\map\n^-$ (resp.
$\CN\hookrightarrow\g\map\n^+$). In particular we obtain a
description of the algebra $\Ext_{\u_\ell}(\underC,\underC)$.

     \sssn \Cor \label{freee}
$\Ext_{\u_{\ell}}(\underC,\underC)=S\bul(\n^{-*})\ten H_0\ten
S\bul(\n^{+*})$ as a module over $S\bul(\n^{-*})\ten
S\bul(\n^{+*})$. Here $H_0:=\bigoplus_{w\in W}\CC\xi_w$. \qed

     Ginzburg and Kumar also proved the following statements
(see \cite{GK}). \begin{itemize} \item[(i)]
$\Ext_{\u_\ell^-}(\underC,\underC)=\underset{w\in
W}{\bigoplus}\CC\nu^+_w\ten S\bul(\n^{-*})$ as a $X\times
\Z$-graded vector space. Here the homological grading of the
element $\nu_w^+$, $w\in W$, is equal to $l(w)$, and its
$X$-grading equals $\rho-w(\rho)$. The homological grading
(resp. the $X$-grading) of the element $f_\alpha^*\in\n^{-*}$
equals $2$ (resp. $\ell\alpha'$). A similar statement holds for
the $\ext$ algebra of $\u_\ell^+$. \item[(ii)]
$\Ext_{\b_\ell^-}(\underC,\underC)$ (resp.
$\Ext_{\b_\ell^+}(\underC,\underC)$) is equal to
$S\bul(\n^{-*})$ (resp. $S\bul(\n^{+*})$) both as an
associative algebra and as a $\n^-$- (resp. $\n^+$-) module.
\item[(iii)] The natural map $\Ext_{\u_\ell}
(\underC,\underC)\map \Ext_{\b_\ell^-}(\underC,\underC)$
(resp. $\Ext_{\u_\ell}(\underC,\underC)\map
\Ext_{\b_\ell^+}(\underC,\underC)$) given by the restriction
functor coincides with the morphism of the algebras of
functions provided by the inclusion of the affine manifolds
$\n^-\hookrightarrow\CN$ (resp. $\n^+\hookrightarrow\CN$).
\end{itemize} \label{A+-} We set $A_+:=\underset{w\in
W}{\bigoplus}\CC\nu^+_w\subset
\Ext_{\u_\ell^-}(\underC,\underC)$ and $A_-:=\underset{w\in
W}{\bigoplus}\CC\nu^-_w\subset
\Ext_{\u_\ell^+}(\underC,\underC)$. Note that the restriction
functor $\Res_{\u_\ell^+}^{\b_\ell^+}$ (resp.
$\Res_{\u_\ell^-}^{\b_\ell^-}$) provides the inclusion of
algebras $$
\Ext_{\b_\ell^+}(\underC,\underC)=S\bul(\n^{+*})\map
\Ext_{\u_\ell^+}(\underC,\underC) \text{ (resp. }
\Ext_{\b_\ell^-}(\underC,\underC)=S\bul(\n^{-*})\map
\Ext_{\u_\ell^-}(\underC,\underC)). $$ \sssn \Lemma
\begin{itemize} \item[(i)] The subspace $A_-$ (resp. $A_+$) is
a subalgebra in $\Ext_{\u_\ell^+}(\underC,\underC)$ (resp. in
$\Ext_{\u_\ell^+}(\underC,\underC)$). \item[(ii)] The
multiplication maps provide the vector space isomorphisms $$
\Ext_{\b_\ell^-}(\underC,\underC)\ten A_+\til{\map}
\Ext_{\u_\ell^-}(\underC,\underC) \text{ and }
\Ext_{\b_\ell^+}(\underC,\underC)\ten A_-\til{\map}
\Ext_{\u_\ell^+}(\underC,\underC).\qed $$ \end{itemize} Note
that the algebra $\u_\ell^0$ acts trivially on the first
factors in these decompositions.

     \sssn \Prop \label{afrob} The algebras $A_-$ and $A_+$ are
Frobenius.

     \dok We consider only the case of $A_+$. Recall that when
calculating cohomology of the algebra $\u_\ell^-$ Ginzburg and
Kumar used certain degenerations of the algebra defined by Kac,
De Concini and Procesi in \cite{DCKP} with the help of quantum
PBW filtration on $\u_\ell^-$. The algebra
$\gr^{PBW}(\u_\ell^-)$ appears to be a certain quotient algebra
of the {\em quantum symmetric algebra} with the set of
generators equal enumarated by the standard basis in the space
$\n^-$. In particular it was proved in \cite{GK} that $\Ext_{\gr^{PBW}
(\u_\ell^-)}(\underC,\underC)=\Lambda\bul_\ell(\n^
{-*})\ten
     S\bul(\n^{-*})$ as an associative algebra. Here
$\Lambda_\ell\bul(\n^{-*})$ denotes the {\em quantum exterior
algebra} on generators enumerated by elements of the standard
basis in $\n^{-*}$. In particular this algebra is Frobenius
with the trace map given by the projection on the top
homological grading component (that is one dimensional).

     Consider now the spectral sequence with the term $E_2$
equal to $\Ext_{\gr^{PBW}(\u_\ell^-)}(\underC,\underC)$ that
converges to $\Ext_{\u_\ell^-}(\underC,\underC)$. Note that for
$\ell$ big enough the subalgebra
$\Lambda_\ell\bul(\n^{-*})\subset
\Ext_{\gr^{PBW}(\u_\ell^-)}(\underC,\underC)$ is a
DG-subalgebra. Moreover its images on the next terms of the
spectral sequence are DG-subalgerbas too. Now note that the
subalgebra $\gr A_+$ in the term $E_\infty$ of the spectral
sequence comes from the described DG-subalgebra in the term
$E_2$. In particular the one dimensional top homological
grading component of $\Lambda_\ell\bul(\n^{-*})$ has its
nonzero representative on the level $E_\infty$. Thus it has its
nonzero representatives on all the levels $E_m$, $m\ge 2$.

     Consider the algebra $\gr^F A_+$ appearing in the
$E_\infty$ term of the spectral sequence. Put the linear
functional generating the trace map on $\gr^F A_+$ equal to the
projection on the one dimensional top homological grading
component. We claim that the corresponding trace pairing on the
algebra is nondegenerate. Indeed, it is so on the algebra
$\Lambda_\ell\bul(\n^{-*})$ appearing in the $E_2$ term of the
spectral sequence. On the other hand note that the $X$-grading
components of the space $\Lambda_\ell\bul(\n^{-*})$ with the
weights $\rho-w(\rho)$, $w\in W$, are one dimensional. Moreover
the subspace $\bigoplus_{w\in
W}\left(\Lambda_\ell\bul(\n^{-*})\right)_{\rho-w(\rho)}$ is
self-dual with respect to the trace pairing on
$\Lambda\bul_\ell(\n^{-*})$. Thus the term $E_\infty$ of the
spectral sequence becomes an orthogonal direct summ and in the
term $E_2$. We have proved that the algebra $\gr^FA_+$ is
Frobenius. Finally by Lemma 2.4.4 from \cite{Ar1} tha algebra
$A_+$ is Frobenius itself. \qed

     \subsection{Cohomology of the algebra $\til{\u}_{\ell}$.}
\label{degen} Note that the algebra $\til{\u}_\ell$ contains
the algebra $\u^-_\ell\ten\u_\ell^+$ {\em as a subalgebra},
moreover, this subalgebra is {\em normal} in $\til{\u}_\ell$
with the quotient algebra equal to $\u_\ell^0$. On the other
hand, by standard arguments, we have $$
\Ext_{\u_\ell^-\ten\u_\ell^+}(\underC,\underC)=
\Ext_{\u_\ell^-}(\underC,\underC)\ten
\Ext_{\u_\ell^+}(\underC,\underC) $$ as an associative algebra.
Now recall that the algebra $\u_\ell^0$ is semisimple being a
group algebra of a finite group. We have proved the following
statement.

     \sssn \Lemma $\Ext_{\til{\u}_\ell}(\underC,\underC)=
\left( \Ext_{\u_\ell^-}(\underC,\underC)\ten
\Ext_{\u_\ell^+}(\underC,\underC) \right)^{\u_\ell^0}$ as a
$\ell\cdot X\times\Z$-graded vector space. Here
$(\cdot)^{\u_\ell^0}$ denotes taking $\u_\ell^0$-invariants.
Moreover, the restriction functor
$\Res_{\u_\ell^-\ten\u_\ell^+}^{\til{\u}_\ell}$ provides an
isomorfism of associative algebras
$\Ext_{\til{\u}_\ell}(\underC,\underC)\cong\left(\Ext_{\u_\ell^
-\ten\u_\ell^+}
     (\underC,\underC)\right)^{\u_\ell^0}$. \qed

     The following proposition is proved similarly to Theorem
3.1 from \cite{Ar5}.

     \sssn \Lemma $ \ch\left( \left(
\Ext_{\u_\ell^-}(\underC,\underC)\ten
\Ext_{\u_\ell^+}(\underC,\underC) \right)^{\u_\ell^0},t\right)=
{\frac{\underset{w\in W}{\sum}t^{2l(w)}} {\underset{\alpha\in
R^+}{\prod}(1-e^{-\ell\alpha}t^2)(1-e^{\ell\alpha}t^2)}}. $
Here the indeterminate $t$ stands for the homogeneous
$\Z$-grading, for a weight $\alpha=\sum a_ii'$ the symbol
$e^\alpha$ denotes the monomial $\prod_i(e^{i'})^{a_i}$. \qed

     Let us denote the filtration on the space
$\Ext_{\u_\ell}(\underC,\underC)$ that corresponds to the
filtration $F$ on the algebra $\u_\ell$ by the same letter.

     \sssn \Cor \label{decomp} $\gr^F\Ext_{\u_\ell}(\underC,\underC)
     = \Ext_{\til{\u}_\ell}(\underC,\underC)$ as a graded
associative algebra.\qed

     Summing up the previous considerations we obtain the
following statement.

     \sssn \Prop \label{cohtil}
$\Ext_{\til{\u}_\ell}(\underC,\underC)=S\bul(\n^{-*})\ten\gr^FH
_0\ten
S\bul(\n^{+*})$ as an associative algebra. Here the imbedding
$S(\n^{-*})\hookrightarrow \Ext_{\til{\u}_\ell}(\underC,\underC)$
resp. the imbedding $S(\n^{+*})\hookrightarrow
\Ext_{\til{\u}_\ell}(\underC,\underC)$) is provided by the projection
of algebras $\til{\u}_\ell\map\b_\ell^-$: $\u_\ell^+\map 0$ (resp.
$\til{\u}_\ell\map\b_\ell^+$: $\u_\ell^-\map 0$).  \qed

Comparing grading components in the previous equality with gradings
far smaller than $\ell$ we obtain the following statement.

\sssn \Cor
The associative algebra $\gr^FH_0=\left(A_-\ten
A_+\right)^{\u_\ell^0}$.\qed

\section{$\protect\Exts_{\u_\ell}(\protect\underC,\protect\underC)$
as a
$\protect\Ext_{\u_\ell}(\protect\underC,\protect\underC)$-module.}
First we recall briefly the general setup for semiinfinite
cohomology of the small quantum group. The definition of
semiinfinite cohomology presented below is a specialization of
the general one in the case of a finite dimensional graded
algebra $\u=\b^-\ten\u^+$ such that $\b^-$ is nonpositively
graded, and $\u^+$ is a positively graded {\em Frobenius}
algebra with $\u^+_0=\CC$.

\subsection{Definition of semiinfinite cohomology.}\label{setup}
Consider first the {\em semiregular} $\u$-bimodule
$S_\u^{\u^+}=\u\ten_{\u^+}\u^{+*}$, with  the right $\u$-module
structure provided by the isomorphism of left $\u^+$-modules
$\u^+\cong\u^{+*}$.  Note that $S_{\u}^{\u^+}$ is free over the
algebra $\u$ both as a right and as a left module.

For a complex of graded
$\u$-modules $M\bul=\underset{p,q\in\Z}
{\bigoplus}M_p^q,\ d:\ M_p^q\map M_p^{q+1}$ we
define the support of $M\bul$ by $\supp
M\bul:=\{(p,q)\in\Z^2|M_p^q\ne0\}$. We say that a complex $M\bul$ is
{\em concave} (resp. {\em convex}) if there exist $s_1,s_2\in\N,
t_1,t_2\in\Z$ such that $\supp M\bul\
subset\{(p,q)\in\Z^2| s_1q+p\le
t_1, s_2q-p\le t_2\}$ (resp. $\supp M\bul\subset\{(p,q)\in\Z^2|
s_1q+p\ge t_1, s_2q-p\ge t_2\}$).

Let $M\bul,N\bul\in{\mathcal C}om(\u\mod)$. Suppose that  $M\bul$ is
convex and $N\bul$ is concave. Choose a convex (resp. concave) complex
$R_{\uparrow}\bul(M\bul)$ (resp.  $R_{\downarrow}\bul(N\bul)$) in
${\mathcal C}om(\u\mod)$  quasiisomorphic to $M\bul$ (resp. $N\bul$)
and consisting of $\u^+$-free (resp.  $\u^-$-free)
modules.

\sssn
\Def
We set
$$
\Exts_{\u}(M\bul,N\bul)
:=H\bul(\Hom_{\u}(R\bul_{\uparrow}(M\bul),S_{\u}^{\u^+}\ten_\u
R_{\downarrow}\bul(N\bul))).
$$
\sssn
\Lemma (see. [Ar1] Lemma 3.4.2, Theorem 5)
The spaces  $\Exts_{\u}(M\bul,N\bul)$ do not depend on the choice of
resolutions and define functors
$$
\ext_{\u}^{\si+k}(\cdot,\cdot):\ \u\mod\times\u\mod\map{\mathcal
V}ect,\ k\in\Z.\qed
$$
Below we consider semiinfinite cohomology of
algebras $\u_\ell$, $\til{\u}_\ell$, $\b^+_\ell$, $\b_\ell^-$ etc.
with coefficients in $X$-graded modules. The $\Z$-grading on such a
module is obtained from the $X$-grading using the function $\hgt:\
X\map \Z$.

Evidently the spaces $\Exts_{\u_{\ell}}(M\bul,N\bul)$
(resp.  $\Exts_{\til{\u}_\ell}(\til{M}\bul,\til{N}\bul)$) posess
natural $\ell\cdot X$-gradings.
The following statement is a direct consequence of Lemma~\ref{mainu}(ii).

\sssn
\Lemma  \label{action}
Let $M,N\in\u_\ell\mod$ be restrictions of some $\U_\ell$-modules.
Then the spaces $\Exts_{\u_{\ell}}(M,N)$ have  natural structures of
$\g$-modules, and the  $\ell\cdot X$-gradings on them coincide with the
$X$-gradings provided by the weight decompositions of the modules with
respect to the standard Cartan subalgebra in $\g$.\qed

In particular we consider the character of this $\g$-module.

\subsection{Semiinfinite cohomology of the trivial $\til{\u}_\ell$-module}
The following statement sums up the main results from \cite{Ar5}.
\vskip 1mm
\noindent
\Theorem                            \label{main}
\begin{itemize}
\item[(i)]
$\Exts_{\til{\u}_\ell}(\underC,\underC)=\left(\Exts
_{\u_\ell^-\ten\u_\ell^+}
(\underC,\underC)\right)^{\u_\ell^0}.$
\item[(ii)]
$
\ch\left( \Exts_{\til{\u}_\ell}(\underC,\underC),t\right)=
t^{-\sharp(R^+)}e^{-2\ell\rho}{\frac{\underset{w\in W}{\sum}t^{2l(w)}}
{\underset{\alpha\in
R^+}{\prod}(1-e^{-\ell\alpha}t^{-2})(1-e^{-\ell\alpha}t^2)}}.  $

The right hand side of the equality is considered
as an element in
$\CC[t,t^{-1}][[e^{-i'},i\in I]]$.
\item[(iii)]
Let $M$ be a filtered module over the filtered algebra $\u_\ell$ with
the filtration $F$.  Then there exists a spectral sequence with the
term $E_1=\Exts_{\til{\u}_\ell}(\underC,\gr M)$, converging to
$E_{\infty}=\gr^F\Exts_{\u_\ell}(\underC,M)$.  \item[(iv)] For the
trivial $\u_\ell$-module the described spectral sequence degenerates
in the term $E_1$. In particular  as a $\ell\cdot X\times \Z$-graded
vector space $\Exts_{\til{\u}_\ell}(\underC,\underC)=
\Exts_{\u_\ell}(\underC,\underC)$.\qed \end{itemize} Let us compare
semiinfinite cohomology of the trivial $\b_\ell^+$-module with the
usual cohomology of $\b_\ell^+$ with coefficients in this module.

\sssn
\Lemma \label{compare}
     $\Exts_{\b_\ell^+}(\underC,\underC)=\Ext_{\b_\ell}(\underC
,\CC((-\ell+1)2\rho))$.

     \dok The triangular decomposition of the algebra
$\b_\ell^+$ is as follows: $\b_\ell^+=\u_\ell^+\ten\u_\ell^0$.
Then by definition of semiinfinite cohomology, using the fact
that $\u_\ell^0$ semisimple, we see that $$
\Exts_{\b^+_\ell}(\underC,\underC)=
H\bul\left(\Hom_{\b^+_\ell}(R\bul_{\uparrow}(\underC),S_{\b^+_\ell}
^{\u_\ell^+}\ten_{\u_\ell}
\underC)\right),
$$
     where $R\bul_{\uparrow}(\underC)$ denotes a concave
$\u_\ell^+$-free resolution of the trivial $\b_\ell^+$-module.
Next note that the $\b^+_\ell$-module
$S_{\b_\ell^+}^{\u_\ell^+}\ten_{\b_\ell^+}\underC=
\b_\ell^{+*}\ten_{\b_\ell^+}\underC=\CC(-(\ell-1)2\rho)$. Thus
we have \begin{gather*} \Exts_{\b^+_\ell}(\underC,\underC)=
H\bul\left(\Hom_{\b^+_\ell}(R\bul_{\uparrow}(\underC),\CC(-(\ell
-1)2\rho))\right)\\
     =\left(H\bul(\Hom_{\u^+_\ell}(R\bul_{\uparrow}(\underC),\C
C(-(\ell-1)2\rho)))\right) ^{\u_\ell^0}
     =\Ext_{\b^+_\ell}(\underC,\CC(-(\ell-1)2\rho)).\qed
\end{gather*}
Next we investigate the structure of the
$\Ext_{\til{\u}_\ell}(\underC,\underC)$-module
$\Exts_{\til{\u}_\ell}(\underC,\underC)$.

\sssn
\Prop    \label{module}
Up to the homological grading shift by $\sharp(R^+)$ we have
$$
\Exts_{\til{\u}_\ell}(\underC,\underC)=
\Coind_{\Ext_{\b_\ell^+}(\underC,\underC)}^{
\Ext_{\til{\u}_\ell}(\underC,\underC)}
\left(\Ext_{\b_\ell^+}(\underC,\CC(-(\ell-1)2\rho))\right)
$$
as a
$\Exts_{\til{\u}_\ell}(\underC,\underC)$-module.

\dok
     By Theorem \ref{main}, Lemma~\ref{compare} and
Corollary~\ref{decomp} we have \begin{gather*}
\Exts_{\til{\u}_\ell}(\underC,\underC)=
\left(\Exts_{\u_\ell^-\ten\u_\ell^+}
(\underC,\underC)\right)^{\u_\ell^0}=
\left(\Tor^{\u_\ell^-}(\underC,\underC)\ten\CC((1-\ell)2\rho)\ten
\Ext_{\u_\ell^+}
     (\underC,\underC)\right)^{\u_\ell^0}\\
=\left(S\bul(\n^-)\ten A_+^*\ten S\bul(\n^{+*})\ten
A_-\ten\CC((1-\ell)2\rho)\right)^{\u_\ell^0}\\ =
S\bul(\n^-)\ten
S\bul(\n^{+*})\ten\left(A_+^*\ten\CC((1-\ell)2\rho)\ten
A_-\right)^{\u_\ell^0}. \end{gather*} Now recall that the
algebra $A_+$ is Frobenius. Making the isomorphism $A_+=A_+^*$
compatible with the $\u_\ell^0$-action, we see that
$A_+=A_+^*\ten\CC(-2\rho)$. Thus we obtain \begin{gather*}
\Exts_{\til{\u}_\ell}(\underC,\underC)= S\bul(\n^-)\ten
S\bul(\n^{+*})\ten
\left(A_+\ten\CC(-2\rho)\ten\CC((1-\ell)2\rho)\ten
A_-\right)^{\u_\ell^0}\\ = S\bul(\n^-)\ten \gr^F H_0\ten
S\bul(\n^{+*})\ten\CC(-2\ell\rho). \end{gather*} Note that the
algebra $\gr H_0$ is Frobenius itself, thus the latter
$\Exts_{\til{\u}_\ell}(\underC,\underC)$-module can be
considered both as $$
\Coind_{\Ext_{\b_\ell^+}(\underC,\underC)}^{
\Ext_{\til{\u}_\ell}(\underC,\underC)}
\left(\Ext_{\b_\ell^+}(\underC,\CC(-2\ell\rho))\right) \text{
and as } \Ind_{\Ext_{\b_\ell^-}(\underC,\underC)}^{
\Ext_{\til{\u}_\ell}(\underC,\underC)}
\left(\Tor^{\b_\ell^-}(\underC,\CC(-2\ell\rho))\right). $$
Finally precise calculation of homological gradings shows that
the grading on the firs t module should be shifted by
$-\sharp(R^+)$, and the grading on the second one should be
shifted by $\sharp(R^+)$. Thus we have \begin{gather*}
\Exts_{\til{\u}_\ell}(\underC,\underC)=
\Coind_{\Ext_{\b_\ell^+}(\underC,\underC)}^{
\Ext_{\til{\u}_\ell}(\underC,\underC)}
\left(\Ext_{\b_\ell^+}(\underC,\CC(-(\ell-1)2\rho))\right)\\ =
\Ind_{\Ext_{\b_\ell^-}(\underC,\underC)}^{
\Ext_{\til{\u}_\ell}(\underC,\underC)}
\left(\Tor^{\b_\ell^-}(\underC,\CC(-(\ell-1)2\rho))\right).\qed
\end{gather*} \sssn \Cor
$\Exts_{\til{\u}_\ell}(\underC,\underC)$ is both free over the
algebra $\Ext_{\b_\ell^+}(\underC,\underC)$ and co-free over
the algebra $\Ext_{\b_\ell^-}(\underC,\underC)$. \qed

     \subsection{Functions on $\CN$ with singularities along
coordinate hyperplanes.} Consider the divisors on the nilpotent
cone of the form $D_\alpha:=\{n\in\CN|f_\alpha^*(n)=0\}$, where
$f_\alpha^*$ denote the standard basic linear functions on
$\n^-$, $\alpha\in R^+$. Note that we have
$\CN\supset\n^+=\underset{\alpha\in R^+}{\bigcap}D_\alpha$.

     Consider the open subset in the nilpotent cone of the form
$ \CN_{(f_{\alpha_1}^*,\ldots,f_{\alpha_{\sharp(R^+)}}^*)}:=
\CN\setminus\underset{\alpha\in R^+}{\bigcup}D_\alpha$ and the
$H^0(\CN,\O_{\CN})$-module $H^0(
\CN_{(f_{\alpha_1}^*,\ldots,f_{\alpha_{\sharp(R^+)}}^*)},\O_{\C
N})$.

     It is well known how to obtain this
$H^0(\CN,\O_{\CN})$-module from the free module
$H^0(\CN,\O_{\CN})$ via an inductive limit construction.

     For $\alpha\in R^+$ consider the
$H^0(\CN,\O_{\CN})$-module map $\mu_{\alpha}:\
H^0(\CN,\O_{\CN})\map H^0(\CN,\O_{\CN})$ given by multiplicaton
by the element $f_\alpha^*\in\n^{-*}$.

     Next we construct an inductive system of
$H^0(\CN,\O_{\CN})$-modules ${\mathcal I}_\bullet$ enumerated
by the partially ordered set $\Z_+^{R^+}$ with ${\mathcal I}_v=
H^0(\CN,\O_{\CN})$ for any $v\in \Z_+^{R^+}$ and with morphisms
$$ {\mathcal I}_{(\ldots,v_\alpha,\ldots)}\map {\mathcal
I}_{(\ldots,v_\alpha+1,\ldots)} $$ equal to $\mu_{\alpha}$.

\sssn
\Lemma
The
$H^0(\CN,\O_{\CN})$-module
$H^0(
\CN_{(
f_{\alpha_1}^*,\ldots,f_{\alpha_{\sharp(R^+)}}^*)},
\O_{\CN})$ is equal to
$\dirlim {\mathcal I}_{\bullet}$.\qed

Note that the
$H^0(\CN,\O_{\CN})$-module
$H^{\sharp R^+}_{\n^+}(\CN,\O_{\CN})$
can be realised naturally as a bottom subquotient module
in
$H^0(
     \CN_{(f_{\alpha_1}^*,\ldots,f_{\alpha_{\sharp(R^+)}}^*)},
\O_{\CN})$ under the filtration on the latter module by
singularities of functions. More precisely, recall that
$H^0(\CN,\O_{\CN}) =S\bul(\n^{-*})\ten H_0\ten
S\bul(\n^{+*})$as a module over $S\bul(\n^{+*})\ten
S\bul(\n^{-*})$. Fix the set of standard generators in the
algebra $S(\n^{-*})$: $S(\n^{-*})=\CC[f_\alpha^*,\ \alpha\in
R^+]$. Then by the previous Lemma $H^0(
\CN_{(f_{\alpha_1}^*,\ldots,f_{\alpha_{\sharp(R^+)}}^*)},
\O_{\CN})$ is equal to $\CC[f_\alpha^{*\pm1},\alpha\in R^+]\ten
H_0\ten S\bul(\n^{+*})$. On the other hand note that $H^{\sharp
R^+}_{\n^+}(\CN,\O_{\CN})$ is isomorphic to
$\underset{\alpha\in
R^+}{\bigotimes}\CC[f_\alpha^{*\pm1}]/\CC[f_\alpha^*]\ten
H_0\ten S\bul(\n^{+*})$.

Our main goal here is to construct an isomorphism of
$H^0(\CN,\O_{\CN})$-modules
$$H^{\sharp R^+}_{\n^+}(\CN,\O_{\CN})\til{\map}
\Exts_{\u_\ell}(\underC,\underC).$$
     We start with constructing a system of
$H^0(\CN,\O_{\CN})$-module morphisms $\varphi_\bullet:\
{\mathcal I}_{\bullet}\map \Exts_{\u_\ell}(\underC,\underC)$
that would provide a morphism from the direct limit of the
inductive system to the semiinfinite cohomology module.

     \subsection{Construction of the morphisms
$\varphi_\bullet$.} Note that since each module ${\mathcal
I}_v$, $v\in \Z_+^{R^+}$, is free over the algebra
$H^0(\CN,\O_{\CN})$ with one generating element, to construct
the system of morphisms from the described inductive system to
$\Exts_{\u_\ell}(\underC,\underC)$ one has to specify the
images of $1\in H^0(\CN,\O_{\CN})$ under the morphisms
$\varphi_v$, $v\in\Z_+^{R^+}$. Then
$\varphi_v(a)=a\cdot\varphi_v(1)$, where $a\in
H^0(\CN,\O_{\CN})=\Ext_{\u_\ell}(\underC,\underC)$ (the latter
algebra acts naturally on $\Exts_{\u_\ell}(\underC ,\underC)$).
Suppose that $$ \varphi_{(\ldots,v_\alpha,\ldots)}:\ 1\mapsto
m_{(\ldots,v_\alpha,\ldots)}\in
\Exts_{\u_\ell}(\underC,\underC);\
\varphi_{(\ldots,v_\alpha+1,\ldots)}:\ 1\mapsto
m_{(\ldots,v_\alpha+1,\ldots)}\in
\Exts_{\u_\ell}(\underC,\underC). $$ Then evidently the only
condition on the sequence of elements
$\{m_v\}_{v\in\Z_+^{R^+}}$ is that $$ f_\alpha^*\cdot
m_{(\ldots,v_\alpha+1,\ldots)}=m_{(\ldots,v_\alpha,\ldots)}. $$
Consider the element $$ \nu_e^+\ten\nu_e^-\in (A_+\ten
A_-)^{\u^0}\ten \CC(-2\ell\rho)=
\gr^FH_0\ten\CC(-2\ell\rho)\subset
\Exts_{\til{\u}_\ell}(\underC,\underC). $$ Here we preserve
notanion from \ref{A+-} for the base elements in the algebras
$A_+$ and $A_-$, and $e\in W$ denotes the unity element.

     Note that the corresponding $\ell\cdot X\times\Z$-grading
component in the space $\Exts_{\til{\u}_\ell}(\underC,\underC)$
is one dimensional. Since the filtration $F$ on the space
$\Exts_{{\u}_\ell}(\underC,\underC)$ is well defined with
respect to the $\ell\cdot X\times\Z$-grading, we can consider
the described vector as an element of
$\ext_{\u_\ell}^{\frac\infty2-\sharp(R^+)}(\underC,\underC)$.

     We are starting to construct the system of morphisms
$\varphi_v$, $v\in \Z_+^{R^+}$. Put $$\varphi_{(0,\ldots,0)}:\
\Ext_{\u_\ell}(\underC,\underC) [\sharp(R^+)] \map
\Exts_{\u_\ell}(\underC,\underC),\
\varphi_{(0,\ldots,0)}(1)=\nu_e^-\ten\nu_e^+. $$ \sssn \Theorem
\label{cofree} $\Exts_{{\u}_\ell}(\underC,\underC)$ is co-free
over the algebra
$S\bul(\n^{-*})=\Ext_{\b_\ell^-}(\underC,\underC)$ with the
space of co-generators equal to $H_0\ten\CC(-2\ell\rho)\ten
S\bul(\n^{+*})$.

     \dok First recall that by~\ref{freee} the
$\Ext_{\b_\ell^+}(\underC,\underC)$-module
$\Exts_{\til{\u}_\ell}(\underC,\underC)$ is free with the space
of generators equal to $$S\bul(\n^-)\ten
H_0\ten\CC(-2\ell\rho).$$ On the other hand $
\Exts_{\til{\u}_\ell}(\underC,\underC)=
\gr^F\Exts_{{\u}_\ell}(\underC,\underC)$, and this equality is
an isomorphism of modules over $\Ext_{\b_\ell^
+}(\underC,\underC) =\gr^F\Ext_{\b_\ell^+}(\underC,\underC)$.
Fix the set of linear independent elements $\{\overline{b}_p|
p\in\Z_+^{R^+}\times W\}$ in the generators space of
$\Exts_{\til{\u}_\ell}(\underC,\underC)$ and choose the set of
representatives of the elements $\{b_p| p\in\Z_+^{R^+}\times
W\}\subset \Exts_{\u_\ell}(\underC,\underC)$. It is easy to
verify that the $\Ext_{\b_\ell^+}(\underC,\underC)$-submodule
in $\Exts_{\u_\ell}(\underC,\underC)$ generated by this set is
free. Moreover its character coincides with the one of
$\Exts_{\u_\ell}(\underC,\underC)$, thus the module of
semiinfinite cohomology is
$\Ext_{\b_\ell^+}(\underC,\underC)$-free itself. Now to
complete the proof of the Theorem recall the following
statement from~\cite{Ar1}.

     Denote by $\hat{\u}_\ell$ the finite quantum group defined
in the same way as $\u_\ell$, but with $\zeta$ replaced by
$\zeta^{-1}$ in the defining relations.

     \sssn \label{pairing} \Lemma (see \cite{Ar1}, Proposition
5.4.3) There exists a nondegenerate contragradient pairing $$
\langle \ ,\ \rangle :\ \Exts_{\u_\ell}(\underC,\underC)\times
\Exts_{\hat{\u}_\ell}(\underC,\underC)\map\CC $$ well defined
with respect to the actions of the algebra $H^0(\CN,\O_{\CN})=
\Ext_{\u_\ell}(\underC,\underC)
=\Ext_{\hat{\u}_\ell}(\underC,\underC)$ and of the Lie algebra
$\g$. \qed

     Using this pairing we obtain the set of co-generators
$\{c_p|p\in\Z^{R^+}_+\}$ for the co-free
$\Ext_{\b^-_\ell}(\underC,\underC)$-module
$\Exts_{\u_\ell}(\underC,\underC)$ enumerated by the base
vectors of the space $H_0\ten\CC(-2\ell\rho)\ten
S\bul(\n^{+*})$. The Theorem is proved. \qed

     Note that we can choose the cogenerators $c_p$ constructed
above to be $X\times\Z$-homogeneous. In particular there is a
unique choice of the homogeneous base vector in the space
$$\left(\ext^{\si-\sharp(R^+)}_{\u_\ell}(\underC,\underC)\right
)_{-2\ell\rho}$$
     with the top filtration factor equal to
$\nu_e^-\ten\nu_e^+\in
\ext^{\si-\sharp(R^+)}_{\til{\u}_\ell}(\underC,\underC)$. Now
consid er the direct sum decomposition of the
$\Ext_{\b^-_\ell}(\underC,\underC)$-module $$
\ext^{\si-\sharp(R^+)+\bullet}_{\til{\u}_\ell}(\underC,\underC)
=
     \underset{p\in \Z_+^{R^+}\times W}{\bigoplus}
S\bul(\n^-)\ten c_p $$ and the direct summand $S\bul(\n^-)\ten
(\til{\nu_e^+\ten \nu_e^-})$ where $\til{\nu_e^+\ten \nu_e^-}$
denotes the element representing $\nu_e^+\ten\nu_e^-$ . Fixing
standard generators $f_\alpha$, $\alpha\in R^+$ in the space
$\n^-$ we obtain a linear base in $S\bul(\n^-)\ten (\nu_e^+\ten
\nu_e^-)$ of the form $m_v:=f_{\alpha_1}^{v_1}\ldots\ldots
f_{\alpha_{\sharp(R^+)}}^{v_{\sharp(R^+)}}\ten(\til{\nu_e^+\ten
\nu_e^-})$.

     \sssn \Cor We have $f_\alpha^*\cdot
m_{(\ldots,v_\alpha+1,\ldots)}=m_{(\ldots,v_\alpha,\ldots)}$.
\qed

     Thus we obtain a sequence of morphisms of
$\Ext_{\u_\ell}(\underC,\underC)$-modules $$ \varphi_v:\
{\mathcal I}_v= \Ext_{\u_\ell}(\underC,\underC)
[\sharp{R^+}+2|v|] \map \Exts_{\u_\ell}(\underC,\underC),\
\varphi_v(1)= m_v,\ v\in\Z_+^{R^+}. $$ Here
$|v|=v_1+\ldots+v_{\sharp(R^+)}$. We have constructed a
morphism $$
H^0(\CN_{(f_{\alpha_1}^*,\ldots,f_{\alpha_{\sharp(R^+)}}^*)},\O
_{\CN})=
     \dirlim{\mathcal I}_\bullet\map
\Exts_{\u_\ell}(\underC,\underC). $$ \subsection{Construction
of the isomorphism
$H^{\sharp(R^+)}_{\n^+}(\CN,\O_{\CN})\til{\map}\Exts_{\u_\ell}(
\underC,\underC)$.}
     Our nearest goal is to show that the morphism $\invlim
(\varphi_\bullet)$ provides a map of
$H^0(\CN,\O_{\CN})$-modules
$H^{\sharp(R^+)}_{\n^+}(\CN,\O_{\CN})\map
\Exts_{\u_\ell}(\underC,\underC)$.

     \sssn \Lemma For any $w\in
H^0(\CN_{(f_{\alpha_1}^*,\ldots,f_{\alpha_{\sharp(R^+)}}^*)},\O
_{\CN})
     $ there exists
$v^0=(v_{\alpha_1}^0,\ldots,v_{\alpha_{\sharp(R^+)}}^0)$ such
that for any $v=v^0+v^1$, $v^1\in\Z_+^{R^+}$, the element
$\invlim(\varphi_\bullet)(f_{\alpha_1}^{*v_1}\cdot\ldots\cdot
f_{\alpha_{\sharp(R^+)}}^{*v_{\sharp(R^+)}}w)=0$.

     \dok The statement follows immediately from the
restrictions on $\ell\cdot X\times\Z$-gradings of the space
$\Exts_{\u_\ell}(\underC,\underC)$. \qed

     Recall that $H_{\n ^+}(\CN,\O_{\CN}) $ is the largest
quotient module of
$H^0(\CN_{(f_{\alpha_1}^*,\ldots,f_{\alpha_{\sharp(R^+)}}^*)},\
O_{\CN})$
     supported on $\n^+\subset\CN$.

     \sssn \Cor The map $\invlim(\varphi_\bullet)$ defines a
morphism $\Phi:\ H_{\n^+}^{\sharp(R^+)}(\CN,\O_{\CN})\map
\Exts_{\u_\ell}(\underC,\underC)$. \qed

     In fact it follows from the previous considerations that
for any finite dimensional graded $\u_\ell$-module $M$ the
quasicoherent sheaf on $\CN$ corresponding to the
$H^0(\CN,\O_{\CN})$-module $\Exts_{\u_\ell}(\underC,M)$ is
supported on $\n^+\subset\CN$.

     Note that the maps $\mu_\alpha$, $\varphi_v$ and
$\invlim(\varphi_\bullet)$ are well defined with respect to the
filtration $F$. The maps $\gr^F\mu_\alpha:
\Ext_{\til{\u}_\ell}(\underC,\underC)\map
\Ext_{\til{\u}_\ell}(\underC,\underC)[2] $ are provided by
multiplication by the elements
$f_\alpha^*\in\Ext_{\b_\ell^-}(\underC,\underC)\subset\Ext_{\til{\u}_\ell}
(\underC,\underC)$.
     Moreover we have $$ \gr^F\left(\dirlim{\mathcal
I}_\bullet\right)= \dirlim\left(\gr^F{\mathcal
I}_\bullet\right)= \underset{\alpha\in
R^+}{\bigotimes}\CC[f_\alpha^{*\pm1}]\ten\gr^F H_0\ten
S\bul(\n^{+*}). $$ Note also that the morphisms
$\gr^F\varphi_v:\ \Ext_{\til{\u}_\ell}(\underC,\underC)
[\sharp(R^+)+2|v|] \map \Exts_{\til{\u}_\ell}(\underC,\underC)
$ are defined by $
\gr^F\varphi_v(1)=\overline{m}_v=f_{\alpha_1}^{v_1}\ldots
f_{\alpha_{\sharp(R^+)}}^{v_{\sharp(R^+)}}(\nu_e^-\ten\nu_e^+).
$ It follows that the map $\gr^F\invlim(\varphi_\bullet)$
coincides with the natural one $$ \underset{\alpha\in
R^+}{\bigotimes}\CC[f_\alpha^{*\pm1}]\ten\gr^F H_0\ten
S\bul(\n^{+*})\map \underset{\alpha\in
R^+}{\bigotimes}\left(\CC[f_\alpha^{*\pm1}]/\CC[f_\alpha^*]\right)
\ten\gr^F H_0\ten S\bul(\n^{+*}).
     $$ In particular the map $\gr^F\invlim(\varphi_\bullet)$
is surjective.

     \sssn \Lemma The map $\invlim (\varphi_\bullet)$ is
surjective itself.

     \dok Note that the morphisms $\invlim(\varphi_\bullet)$
and $\gr^F\invlim(\varphi_\bullet)$ preserve $\ell\cdot
X\times\Z$-grading components. The $\ell\cdot
X\times\Z$-grading components of the spaces $\dirlim{\mathcal
I}_\bullet$, $\Exts_{\u_\ell}(\underC,\underC)$,
$\dirlim\left(\gr^F{\mathcal I}_\bullet\right)$ and
$\Exts_{\til{\u}_\ell}(\underC,\underC)$ are finite
dimensional.

     Thus we can check injectivity of the map dual to
$\invlim(\varphi_\bullet)$ instead of surjectivity of
$\invlim(\varphi_\bullet)$. But the injectivity of a filtered
map follows from the injectivity of the associated graded map.
\qed

     Thus we obtain a surjective map $\Phi:\
H^{\sharp(R^+)}_{\n^+}(\CN,\O_{\CN})\map
\Exts_{\u_\ell}(\underC,\underC)$. Next recall the following
statement from \cite{Ar1}.

     \sssn \Prop (see \cite{Ar1}, Theorem A.2.2) $$
\ch\left(H^{\sharp(R^+)}_{\n^+}(\CN,\O_{\CN}),t^2\right)
=\ch\left( \Exts_{\u_\ell}(\underC,\underC),t\right). $$ Here
the indeterminate $t$ in the left hand side of the equality
stands for the grading by the natural action of $\CC^*$, in the
right hand side of the equality it stands for the homological
grading. \qed

     We have proved the Feigin conjectire on the level of
$H^0(\CN,\O_{\CN})$-modules.

     \sssn \Theorem The morphism $\Phi$ constructed above
provides an isomorphism of the $H^0(\CN,\O_{\CN})$-modules $$
\Phi:\ H_{\n^+}^{\sharp(R^+)}(\CN,\O_{\CN})\map
\Exts_{\u_\ell}(\underC,\underC).\qed $$ \section{Isomorphism
of $\g$-modules} Recall that both sides of the
$H^0(\CN,\O_{\CN})$-module isomorphism $\Phi$ carry natural
structures of $\g$-modules.

     The one on $\Exts_{\u_\ell}(\underC,\underC)$ is a
consequence of the fact that the algebra $\u_\ell$ is a normal
subalgebra in $\U_\ell$ with the quotient algebra
$\U_\ell//\u_\ell$ equal to $U(\g)$.

     The $\g$-module structure on
$H_{\n^+}^{\sharp(R^+)}(\CN,\O_{\CN})$ comes from the natural
Lie algebra inclusion $\g\hookrightarrow{\mathcal V}ect(\CN)$
and the standard fact due to Kempf (see \cite{K}) that the Lie
algebra of algebraic vector fields acts on local cohomology
spaces. Below we prove that the two described $\g$-module
structures are isomorphic.

     \subsection{Isomorphism of the $\n^+$-modules.} Consider
the standard positive nilpotent Lie subalgebra $\n^+\subset\g$.
We prove first that $H_{\n^+}^{\sharp(R^+)}(\CN,\O_{\CN})$ and
$\Exts_{\u_\ell}(\underC,\underC)$ are isomorphic as
$\n^+$-modules.

     Note first that the $\n^+$- and
$H^0(\CN,\O_{\CN})$-actions on both
$H_{\n^+}^{\sharp(R^+)}(\CN,\O_{\CN})$ and
$\Exts_{\u_\ell}(\underC,\underC)$ satisfy the equation $$
n\cdot (a\cdot c)= a\cdot(n\cdot c)+[n,a]\cdot c,\ $$ where
$n\in\n^+, a,[n,a]\in H^0(\CN,\O_{\CN})$, $c\in
H_{\n^+}^{\sharp(R^+)}(\CN,\O_{\CN})$ (resp. $c\in
\Exts_{\u_\ell}(\underC,\underC)$) and $[*,*]$ denotes the
natural $\g$-action on $H^0(\CN,\O_{\CN})$.

     Consider the $S\bul(\n^{-*})$-module direct sum
decomposition $$
\ext^{\si-\sharp(R^+)+\bullet}_{{\u}_\ell}(\underC,\underC)=
\underset{p\in \Z_+^{R^+}\times W}{\bigoplus} S\bul(\n^-)\ten
c_p $$ with $c_0=\til{\nu_e^+\ten\nu_e^-}\in
\ext^{\si-\sharp(R^+)}_{{\u}_\ell}(\underC,\underC)$.

     \sssn \Lemma $S\bul(\n^-)\ten c_0\subset
\ext^{\si-\sharp(R^+)}_{{\u}_\ell}(\underC,\underC)$ is a
$\n^+$-submodule.

     \dok We prove the statement of the Lemma by induction by
homological grading. First $\left(S\bul(\n^-)\ten
c_0\right)^{-\sharp(R+)}=\CC c_0$, and the statement follows
from the $X$-grading restrictions on the $\n^+$-action.

     Suppose the statement is proved for $S^k(\n^-)\ten c_0$,
$k<k_0$. We prove that $S^{k_0}(\n^-)\ten c_0$ is a
$\n^+$-submodule in
$\ext^{\si-\sharp(R^+)-2k_0}_{{\u}_\ell}(\underC,\underC)$.

     Let $b\in S^{k_0}(\n^-)$ and $n\in\n^+$. Consider the
element $n\cdot b\in
\ext^{\si-\sharp(R^+)-2k_0}_{{\u}_\ell}(\underC,\underC)$ and
its decomposition $n\cdot b= (n\cdot b)^{\operatorname{yes}}+
(n\cdot b)^{\operatorname{no}}$, where $(n\cdot
b)^{\operatorname{yes}}\in S\bul(\n^-)\ten c_0$ and $ (n\cdot
b)^{\operatorname{no}}\in \underset{p\in \Z_+^{R^+}\times
W,p\ne0}{\bigoplus} S\bul(\n^-)\ten c_p. $ Then for any $\alpha
\in R^+$ the element $f_\alpha^*\cdot(n\cdot
b)=n\cdot(f_\alpha^*\cdot b)+[n,f_\alpha^*]\cdot b$ belongs to
$S\bul(\n^-)\ten c_0$, s ince by induction hypothesis both
summands lie there. Thus for any $\alpha \in R^+$ we have
$f_\alpha^*\cdot (n\cdot b)^{\operatorname{no}}=0$, i.~e. the
element belongs to the $S\bul(\n^{-*})$-invariants subspace of
$\Exts_{{\u}_\ell}(\underC,\underC)$. But by
Theorem~\ref{cofree} the invariants space is equal to $H_0\ten
S\bul(\n^{+*})\ten\CC(-2\ell\rho)$. The latter subspace is
suituated in homological gradings greater than or equal to
$-\sharp(R^+)$. On the other hand the $\n^+$-submodule in
question is situated in homological gradings less than or equal
to $-\sharp(R^+)$. It follows that $ (n\cdot
b)^{\operatorname{no}}=0$. \qed

     \sssn \Cor The action of $\n^+$ on $S\bul(\n^-)\ten c_0$
coincides with the one on $S\bul\left(\g/\b^+\right)$.

     \dok Identify $S\bul(\n^-)\ten c_0$ with
$\hom_\CC(S\bul(\n^{-*}),\CC)$. For such a linear function $b$
note that $b(a)=(a\cdot b)(1)$. Note also that
$S\bul(\n^{-*})=S\bul\left((\g/\b^+)^*\right)$ is a
$\n^+$-submodule in $H^0(\CN,\O_{\CN})$. Thus for $n\in\n^+$,
$b\in \hom_\CC(S\bul(\n^{-*}),\CC)$ we have $$ (n\cdot
b)(a)=a\cdot(n\cdot b)(1)=n\cdot(a\cdot b)(1)-([n,a]\cdot b)(1)
=-b([n,a]).\qed $$

     Recall that the $H^0(\CN,\O_{\CN})$-module
$\Exts_{\u_\ell}(\underC,\underC)$ is generated by the subspace
$S\bul(\n^-)\ten c_0$ and is in fact isomorphic to
$\Ind_{S\bul(\n^{-*})}^{H^0(\CN,\O_{\CN})}S\bul(\n^-)\ten c_0$.
Summing up the previous considerations we note that we have
proved the following statement.

     \sssn \Theorem \label{nplus} Up to an isomorphism there
exists only one $\n^+$-module structure on the space $$
\Exts_{\u_\ell}(\underC,\underC)=H^{\sharp(R^+)}_{\n^+}(\CN,\O_
{\CN})
     $$ well defined with respect to the
$H^0(\CN,\O_{\CN})$-structure. \qed

     \subsection{Tilting modules.} It remains to extend the
described $\n^+$-module isomorphism to a $\g$-module
isomorphism.

     Recall that a finitely generated $\n^+$-locally finite
$\g$-module $M$ diagonazible over $\h$ is called a {\em tilting
module} if it is filtered both by Verma modules and by
contragradient
     Verma modules. H.~ H.~Andersen has proved the following
statement.

     \sssn \Theorem (see \cite{A}) \label{andersen} A tilting
$\g$-module is uniquely determined up to an isomorphism by its
character with respect to the root space decomposition for the
Cartan action. \qed

     Our goal here is to prove that both
$H^{\sharp(R^+)}_{\n^+}(\CN,\O_{\CN})$ and
$\Exts_{\u_\ell}(\underC,\underC)$ are tilting $\g$-modules.

     \subsection{A nondegenerate pairing on
$H^{\sharp(R^+)}_{\n^+}(\CN,\O_{\CN})$} The following
construction was proposed by V.~Ostrik. Instead of constructing
a contragradient pairing on
$H^{\sharp(R^+)}_{\n^+}(\CN,\O_{\CN})$ we obtain a bilinear
pairing $$ H^{\sharp(R^+)}_{\n^+}(\CN,\O_{\CN})\times
H^{\sharp(R^+)}_{\n^-}(\CN,\O_{\CN})\map\CC $$ well defined
with respect to the $\g$-module structures. Consider the
canonical map for the local cohomology spaces provided by the
cup product $$ \langle\ ,\ \rangle:\
H^{\sharp(R^+)}_{\n^+}(\CN,\O_{\CN})\times
H^{\sharp(R^+)}_{\n^-}(\CN,\O_{\CN})\map
H^{2\sharp(R^+)}_{\n^-\cap\n^+}(\CN,\O_{\CN}). $$ \sssn \Lemma
\begin{itemize} \item[(i)] $
H^{2\sharp(R^+)}_{\n^-\cap\n^+}(\CN,\O_{\CN}) =
H^{2\sharp(R^+)}_{\mathbf 0}(\CN,\O_{\CN})=\CC$. Here $\mathbf
0$ denotes the vertex of the nilpotent cone. \item[(ii)] The
map $\langle\ ,\ \rangle$ provides a nondegenerate pairing.\qed
\end{itemize}

     \subsubsection{Springer-Grothendieck resolution of the
nilpotent cone.} To prove that
$H^{\sharp(R^+)}_{\n^+}(\CN,\O_{\CN})$ is $\n^-$-free recall
that in \cite{Ar1} we obtained another geometric realization of
this $\g$-module as follows.

     Consider the simply connected Lie group $G$ with the Lie
algebra equal to $\g$. Choose a maximal torus $H\subset G$
providing the root decomposition of $\g$ and in particular its
triangular decomposition $\g=\n^-\oplus\h\oplus\n^+$. Consider
the Borel subgroup $B\subset G$ with the Lie algebra
$\b^+=\h\oplus\n^+$ and the flag variety $G/B$. The group $G$
acts on $G/B$ by left translations and the restric tion of this
action to $B$ is known to have finitely many orbits. These
orbits are isomorphic to affine spaces and called {\em the
Schubert cells}. The Bruhat decomposition of $G$ shows that the
Schubert cells are enumerated by the Weyl group. Denote the
orbit corresponding to the element $w\in W$ by $S_w$.

     It is well known that the cotangent bundle $T^*(G/B)$ has
a nice realization $ T^*(G/B)=\{(B_x,n)|n\in
\operatorname{Lie}(B_x)\} $, where $B_x$ denotes some Borel
subgroup in $G$ and $n$ is a nilpotent element in the Lie
algebra $\operatorname{Lie}(B_x)$. The map $$ \mu:\
T^*(G/B)\map \CN,\ (B_x,n)\mapsto n, $$ is known to be a
resolution of singularities of $\CN$ called {\em the
Springer-Grothendieck resolution}.

     Recall the following statement from \cite{Ar1}.

     \sssn \Prop (see e. g. \cite{CG}, 3.1.36) \begin{itemize}
\item[(i)] $\mu^{-1}(\n^+)=\underset{w\in
W}{\bigsqcup}T^*_{S_w}(G/B)$, where $T_{S_w}^*(G/B)$ denotes
the conormal bundle to $S_w$ in $G/B$. \item[(ii)]
$H^{\sharp(R^+)}_{\n^+}(\CN,\O_{\CN})\til{\map}
H^{\sharp(R^+)}_{\mu^{-1}(\n^+)}(T^*(G/B),\O_{T^*(G/B)})$ as a
$\g$-module. \qed \end{itemize} \sssn \Cor The $\g$-module
$H^{\sharp(R^+)}_{\n^+}(\CN,\O_{\CN})$ is free over the algebra
$U(\n^-)$. \qed

     In particular it is filtered by Verma modules. Using the
fact that the $\g$-module is self-dual with respect to the
described contragradient pairing we construct a filtration by
contragradient Verma modules on
$H^{\sharp(R^+)}_{\n^+}(\CN,\O_{\CN})$. We have proved the
following statement.

     \sssn \Prop \label{til1} The $\g$-module
$H^{\sharp(R^+)}_{\n^+}(\CN,\O_{\CN})$ is a tilting module.

     \ssn On the other hand recall that in Lemma~\ref{pairing}
we obtained a nondegenerate $\g$-equivariant contragradient
pairing $$ \Exts_{\u_\ell}(\underC,\underC)\times
\Exts_{\hat{\u}_\ell}(\underC,\underC)\map\CC. $$ It follows
from Proposition~\ref{til1} and Theorem~\ref{nplus} that the
$\g$-module $\Exts_{\u_\ell}(\underC,\underC)$ is cofree over
the algebra $U(\n^+ )$. Thus it is filtered by contragradient
Verma modules. Again using the contragradient pairing on the
$\g$-module we produce on it a filtration with subquotients
equal to direct sums of Verma modules.

     We have proved the following statement.

     \sssn \Prop \label{til2} The $\g$-module
$\Exts_{\u_\ell}(\underC,\underC)$ is a tilting module. \qed

     Finally we come to the main statement of the present
paper.

     \sssn \Theorem The $\g$-modules $
\Exts_{\u_\ell}(\underC,\underC)$ and
$H^{\sharp(R^+)}_{\n^+}(\CN,\O_{\CN})$ are isomorphic.

     \dok The statement of the Theorem follows from
Theorem~\ref{andersen} and Propositions~\ref{til1}
and~\ref{til2}. \qed

The Feigin conjecture is proved.

\section{Semiinfinite cohomology of contragradient Weyl modules.}
\ssn
     Fix a dominant integral weight $\lambda\in X$. Consider
the module over the ``big'' quantum group $\U_\ell$ given by $$
\D
W(\lambda):=\left(\Coind_{\B^+_\ell}^{\U_\ell}\CC(\lambda)\right)
^{\operatorname{fin}}
     \text{ (resp. by }
W(\lambda):=\left(\Ind_{\B^+_\ell}^{\U_\ell}\CC(\lambda)\right)
_{\operatorname{fin}}).
     $$ Here $\B_\ell^+=\U_\ell^+\ten\U_\ell^0$ (resp.
$\B_\ell^-=\U_\ell^-\ten\U_\ell^0$ ) denotes the quantum
positive (resp. negative) Borel subalgebra in $\U_\ell$ and
$(*)^{\operatorname{fin}}$ (resp. $(*)_{\operatorname{fin}}$)
denotes the maximal finite dimensional submodule (resp.
quotient module) in $(*)$. The module $\D W(\lambda)$ (resp.
$W(\lambda)$) is called {\em the contragradient Weyl module}
(resp. {\em the Weyl module}) over $\U_\ell$ with the highedst
weight $\lambda$.

     It is known that both modules $W(\lambda)$ and $\D
W(\lambda)$ provide the natural specializations of the finite
dimensional simple modules $L(\lambda)$ over the quantum group
$\U$ at generic values of the quantizing parameter into the
root of unity $\zeta$. In particular we have
$$\ch(W(\lambda))=\ch(\D W(\lambda))=\sum_{w\in
W}{\frac{e^{w\cdot\lambda}}{\prod_{\alpha\in
R^+}(1-e^\alpha)}}, $$ just like in the Weyl character
     formula in the semisimple Lie algebra case. Note also that
$W(0)=\D W(0)=\underC$.

     Below we consider semiinfinite cohomology of the algebra
$\u_\ell$ with coefficients in the contragradient Weyl module
with a $\ell$-divisible highest weight $\ell\lambda$. Ou
considerations were motivated by the following result of
Ginzburg and Kumar (see \cite{GK}).

     Let $p$ denote the projection $T^*(G/B)\map G/B$. Consider
the linear bundle $\L(\lambda)$ on $G/B$ with the first Chern
class equal to $\lambda\in X=H^2(G/B,\Z)$.

     \sssn \Theorem \label{Weyl} \begin{itemize} \item[(i)]
$\ext_{\u_\ell}^{\operatorname{odd}}(\underC,\D
W(\ell\lambda))=0$; \item[(ii)]
$\ext_{\u_\ell}^{2\bullet}(\underC,\D W(\ell\lambda))=
H^0(\CN,\mu_*p^*\L(\lambda))$
     as a $H^0(\CN,\O_{\CN})$-module. \qed \end{itemize}

     The following conjecture provides a natural semiinfinite
analogue for Theorem~\ref{Weyl}.

     \sssn \label{maincon} \Con
     $\ext_{\u_\ell}^{\si+\bullet}(\underC,\D W(\ell\lambda))=
H^0_{\n^+}(\CN,\mu_*p^*\L(\lambda))$
     as a $H^0(\CN,\O_{\CN})$-module. The homological grading
on the left hand side of the equality corresponds to the
grading by the natural action of the group $\CC^*$ on the right
hand side.
     \qed

     \sssn \Cor \label{chformula} \begin{gather*} \ch\left(
\ext_{\u_\ell}^{\si+\bullet}(\underC,\D
W(\ell\lambda)),t\right)\\=
\frac{t^{-\sharp(R^+)}}{\prod_{\alpha\in
R^+}(1-e^{-\ell\alpha})}\sum_{w\in W}
\frac{e^{w(\ell\lambda)}t^{2l(w)}}{\prod_{\alpha\in
R^+,w(\alpha)\in R^+}(1-t^2e^{-\ell\alpha}) \prod_{\alpha\in
R^+,w(\alpha)\in R^-}(1-t^{-2}e^{-\ell\alpha})}.\qed
\end{gather*} Below we present the main steps for the proof of
the conjecture. Details of the proof will be given in the
fourthcoming paper \cite{Ar7}.

     \subsection{Local cohomology with coefficients in $\mu_*
p^*\L(\lambda)$.} For a $H^0(\CN,\O_{\CN})$-module $\M$
consider the natural pairing $ H^i_{\n^+}(\CN,\O_{\CN})\times
H^0(\CN,\M) \map H^i_{\n^+}(\CN,\M). $ Evidently it is
equivariant with respect to the $H^0(\CN,\O_{\CN})$-action.
Thus we obtai n a $H^0(\CN,\O_{\CN})$-module morphism $ s:\
H_{\n^+}^i(\CN,\O_{\CN})\ten_{H^0(\CN,\O_{\CN})}
H^0(\CN,\M)\map H^i_{\n^+}(\CN,\M). $

     \sssn \Prop \label{step1} For $\M=\mu_*p^*\L(\lambda)$ the
map $s$ is an isomorphism. \qed

     \subsubsection{Similar construction for semiinfinite
cohomology.} We will need some more homological algebra. Fix a
graded algebra $A$ with a subalgebra $B\subset A$. Recall that
in \cite{V} and \cite{Ar1} the notion of a complex of graded
$A$-modules K-semijective with respect to the subalgebra $B$
was developed. The following statement gives an analogue of the
standard technique of projective resolutions in the
semiinfinite case.

     \sssn \Theorem (see \cite{Ar1}, Appendix B) Let
$SS_{\u_\ell^+}\bul(*)$ (resp. $SS_{\u_\ell^-}\bul(*)$) denote
a K-semijective concave (resp. convex) resolution of the
$\u_\ell$-module $(*)$ with respect to the subalgebra
$\u_\ell^+$ (resp. $\u_\ell^-$). Then for a finite dimensional
graded $\u_\ell$-module $M$ we have \begin{itemize} \item[(i)]
$H\bul(\Hom_{\u_\ell}( SS\bul_{\u_\ell^-}(\underC),
SS\bul_{\u_\ell^-}(M))=\Ext_{\u_\ell}(\underC,M)$; \item[(ii)]
$H\bul(\Hom_{\u_\ell}( SS\bul_{\u_\ell^+}(\underC),
SS\bul_{\u_\ell^-}(\underC))=\Exts_{\u_\ell}(\underC,M)$. \qed
\end{itemize} \sssn \Cor The composition of morphisms provides
a natural pairing $$
\ext_{\u_\ell}^{\si+i}(\underC,\underC)\times\ext^j_{\u_\ell}(\
underC,M)\map
     \ext^{\si+i+j}_{\u_\ell}(\underC,M).\qed $$ In particular
we obtain a $\Ext_{\u_\ell}(\underC,\underC)$-module map $$
\Exts_{\u_\ell}(\underC,\underC)\ten_{\Ext_{\u_\ell}(\underC,\u
nderC)}\Ext_{\u_\ell}(\underC,M)\map
     \Exts_{\u_\ell}(\underC,M). $$ Combining this construction
with the previous considerations we obtain the following
statement.

     \sssn \Prop There exists a natural
$H^0(\CN,\O_{\CN})$-module morphism $$ \sigma:\
H_{\n^+}^{\sharp(R^+)}(\CN,\mu_*p^*\L(\lambda))\map
\Exts_{\u_\ell}(\underC,\D W(\ell\lambda)). $$ \dok Follows
from Proposition~\ref{step1}. \qed

     Below we show that the morphism $\sigma$ is an
isomorphism. The main
     tool for the demonstration of this fact is the {\em
quasi-BGG complex} providing a specialization of the classical
BGG resolution for a finite dimensional simple $\U$-module
$L(\ell\lambda)$ into the root of unity $\zeta$. This complex
constructed below consists of direct sums of $\U_\ell$-modules
called {\em the quasi-Verma modules}. Moreover its zero
cohomology module equals $W(\ell\lambda)$. On the other hand we
show that semiinfinite cohomology with coefficients in
quasi-Verma modules has a nice geometrical interpretation.

     Now we turn to the construction of the quasi-BGG complex.

     \subsection{Twisted quantum parabolic subalgebras in
$\U_\Q$} Recall that Lusztig has constructed an action of the
{\em braid group} $\CB$ corresponding to the Cartan data
$(I,\cdot)$ by automorphisms of the quantum group $\U_\Q$ well
defined with respect to the $X$-gradings. Fix a reduced
expression of the maximal length element $w_0\in W$ via the
simple reflection elements: $$ w_0=s_{i_1}\ldots
s_{i_{\sharp(R^+)}},\ i_k\in I. $$ Then it is known that this
reduced expression provides reduced expressions for all the
elements $w\in W$: $ w=s_{i_1^w}\ldots s_{i_{l(w)}^w},\ i_k\in
I. $

     Consider the standard generators $\{T_i\}_{i\in I}$ in the
braid group $\CB$. Lifting the reduced expressions for the
elements $w$ from $W$ into $\CB$ we obtain the set of elements
in the braid group of the form $T_w:= T_{i_1^w}\ldots
T_{i_{l(w)}^w}$.

     In particular we obtain the set of {\em twisted Borel
subalgebras} $w(\B_\Q^+)=T_w(\B_\Q^+)\subset\U_\Q$, where
$\B_\Q^+:=\U^+_\Q\ten\U^0_\Q$. Note that
$w_0(\B_\Q^+)=\B_\Q^-=\U_\Q^-\ten\U_\Q^0$.

     Fix a subset $J\subset I$ and consider the {\em quantum
parabolic subalgebra} $\P_{J,\Q}\subset\U_\Q$. By definition
this subalgebra in $\U_\Q$ is generated over $\U_\Q^0$ by the
elements $E_i$, $i\in I$, $F_j$, $j\in J$, and by their quantum
divided powers. The previous construction provides the set of
{\em twisted quantum parabolic subalgebras}
$w(\P_{J,\Q}):=T_w(\P_{J,\Q})$ of the type $J$ with the t wists
$w\in W$.

     Note that the triangular decomposition of the algebra
$\U_\Q$ provides the ones for the algebras $w(\B_\Q^+)$ and
$w(\P_{J,\Q})$: $$ w(\B_\Q^+) = (w(\B_\Q^+))^- \ten \U_\Q^0
\ten (w(\B_\Q^+))^+ \text{ and } w(\P_{J,\Q}) =
(w(\P_{J,\Q}))^- \ten \U_\Q^0 \ten (w(\P_{J,\Q})^+, $$ where
$\left(w(\B_\Q^+)\right)^+ = w(\B_\Q^+)\cap \U_\Q^-$,
$\left(w(\P_{J,\Q})\right)^+ = w(\P_{J,\Q})\cap \U_\Q^-$ etc.

     Specializing the quantizing parameter into the root of
unity $\zeta$ we obtain in particular the subalgebras
$w(\B_\ell^+)\subset\U_\ell$, $w(\P_{J,\ell})\subset \U_\ell$,
$w(\b_\ell^+)\subset\u_\ell$, $w(\p_{J,\ell})\subset \u_\ell$
with the induced triangular decompositions.

     \subsection{Semiinfinite induction and coinduction} From
now on we will use freely the technique of associative algebra
semiinfinite homology and cohomology for a graded associative
algebra $A$ with two subalgebras $B,N\subset A$ equipped with a
triangular decomposition $A=B\ten N$ on the level of graded
vector spaces. We will not recall the construction of these
functors referring the reader to \cite{Ar1} and \cite{Ar2}.

     Let us mention only that these functors are bifunctors
${\mathsf D}(A\mod)\times{\mathsf D}(\oppA\mod)\map{\mathsf
D}({\mathcal V}ect)$ where the associative algebra $\oppA$ is
defined as follows.

     Consider the semiregular $A$-module $S_A^N:=A\ten_NN^*$.
It is proved in \cite{Ar2} that under very weak conditions on
the algebra $A$ the module $S_A^N$ is isomorphic to the
$A$-module $\left(S_A^N\right)':=\hom_B(A,B)$. Thus
$\End_A(S_A^N)\supset N^\opp$ and $\End_A(S_A^N)\supset B^\opp$
as subalgebras. The algebra $\oppA$ is defined as the
subalgebra in $\End_A(S_A^N)$ generated by $B^\opp$ and
$N^\opp$. It is proved in \cite{Ar2} that the algebra $\oppA$
has a triangular decomposition $\oppA=N^\opp\ten B^\opp$ on the
level of graded vector spaces. Yet for an arbitrary algebra $A$
the algebras $\oppA$ and $A^\opp$ do not coincide.

     However the following statement shows that in the case of
quantum groups that corre spond to the root data $(Y,X,\ldots)$
of the {\em finite} type $(I,\cdot)$. the equality of $A^\opp$
and $\oppA$ holds.

     \sssn \Prop We have \begin{itemize} \item[(i)]
$\U^{\sharp}=\U^\opp$, $\U_\Q^{\sharp}=\U_\Q^\opp$,
$\U_\Q^{\sharp}=\U_\Q^\opp$; \item[(ii)] $w(\B_\Q^+)^\sharp=
w(\B_\Q^+)^\opp$, $w(\B_\ell^+)^\sharp= w(\B_\ell^+)^\opp$,
$w(\P_{J,\Q})^\sharp= w(\P_{J,\Q})^\opp$,
$w(\P_{J,\ell})^\sharp= w(\P_{J,\ell})^\opp$. \end{itemize}
\dok The first part is proved similarly to Lemma 9.4.1 from
\cite{Ar6}. The second one follows immediately from the first
one. \qed

     \sssn \Def Let $M\bul$ be a convex complex of
$w(\B_\Q^+)$-modules. By definition set \begin{gather*}
\SInd_{w(\B_\Q^+)}^{\U_\Q}(M\bul):=
\tor_{\si+0}^{w(\B_\Q^+)}(S_{\U_\Q}^{\U_\Q^+},M\bul) \text{ and
} \SCoind_{w(\B_\Q^+)}^{\U_\Q}(M\bul):=
\ext^{\si+0}_{w(\B_\Q^+)}(S_{\U_\Q}^{\U_\Q^-},M\bul).
\end{gather*} The functors
$\SInd_{w(\B_\ell^+)}^{\U_\ell}(\cdot)$,
$\SCoind_{w(\B_\ell^+)}^{\U_\ell}(\cdot)$,
$\SInd_{w(\P_{J,\ell}^+)}^{\U_\ell}(\cdot)$,
$\SCoind_{w(\P_{J,\ell}^+)}^{\U_\ell}(\cdot)$ etc. are defined
in a similar way.

     \sssn \Lemma (see \cite{Ar4}) \begin{itemize} \item[(i)]
$\tor_{\si+k}^{w(\B_\Q^+)}(S_{\U_\Q}^{\U_\Q^+},\cdot)=0$ for
$k\ne0$; \item[(ii)]
$\ext^{\si+k}_{w(\B_\Q^+)}(S_{\U_\Q}^{\U_\Q^-},\cdot)=0$ for
$k\ne0$; \item[(iii)] $\SInd_{w(\B_\Q^+)}^{\U_\Q}(\cdot)$ and
$\SCoind_{w(\B_\Q^+)}^{\U_\Q}(\cdot)$ define {\em exact}
functors $w(\B_\Q^+)\mod\map\U_\Q\mod$.\qed \end{itemize}
Similar statements hold for the algebras $w(\B_\ell^+)$,
$w(\P_{J,\Q})$ and $w(\P_{J,\ell})$.

     \subsection{Quasi-Verma modules} We define the {\em
quasi-Verma module} over the algebra $\U_\Q$ (resp. $\U_\ell$)
with the highest weight $w\cdot\lambda$ by $$ M_\Q^w(w\cdot
\lambda):= \SInd_{w(\B_\Q^+)}^{\U_\Q}(\CC(\lambda)) \text{
(resp. } M_\ell^w(w\cdot \lambda):=
\SInd_{w(\B_\ell^+)}^{\U_\ell}(\CC(\lambda))). $$ The {\em
contragradient quasi-Verma module} $ \D M_\Q^w(w\cdot\lambda)$
resp. $ \D M_\ell^w(w\cdot\lambda)$) is defined by $$ \D
M_\Q^w(\lambda):= \SCoind_{w(\B_\Q^+)}^{\U_\Q}(\CC(\lambda))
\text{ (resp. } \D M_\ell^w(\lambda):=
\SCoind_{w(\B_\ell^+)}^{\U_\ell}(\CC(\lambda))). $$ We list the
main properties of quasi-Verma modules.

     \sssn \Prop (see \cite{Ar4}) \begin{itemize} \item[(i)]
Fix a dominant integral weight $\lambda\in X$. Suppose that
$\xi\in\CC^*$ is not a root of unity. Then the $\U_\xi$-module
$M_\xi^w(w\cdot\lambda):=M_\Q^w(w\cdot\lambda)\ten_{\Q[v,v^{-1}
]}\CC$
     (resp. $\D M_\xi^w(w\cdot\lambda):=\D
M_\Q^w(w\cdot\lambda)\ten_{\Q[v,v^{-1}]}\CC$) is isomorphic to
the usual Verma module $M_\xi(w\cdot\lambda)$ (resp. to the
usual contragradient Verma module $\D M_\xi(w\cdot\lambda)$).
\item[(ii)] For any $\lambda\in X$ we have \begin{gather*}
\ch(M_\Q^w(w\cdot\lambda)) = \ch(\D M_\Q^w(w\cdot\lambda)) =
\ch(M_\ell^w(w\cdot\lambda)) = \ch(\D M_\ell^w(w\cdot\lambda))
= \frac{e^{w\cdot\lambda}}{\prod_{\alpha\in
R^+}(1-e^{-\alpha})}.\qed \end{gather*} \end{itemize} Thus for
a dominant weight $\lambda$ one can consider
$M_\Q^w(w\cdot\lambda)$ as a
     flat family of modules over the quantum group for various
values of the quantizing parameter with the fiber at a generic
$\xi\in\CC^*$ equal to the Verma module $M_\xi(w\cdot\lambda)$.
\vskip 1mm \noindent \Rem Note that by definition
$M_\ell^e(\lambda)=M_\ell(\lambda)$ and $\D
M_\ell^e(\lambda)=\D M_\ell(\lambda)$, where $e$ denotes the
unity element of the Weyl group. In particular for a dominant
weight $\lambda$
     we have a natural projection $M_\ell^e(\lambda)\map
W(\lambda)$ and a natural inclusion $\D
W(\lambda)\hookrightarrow \D M_\ell^e(\lambda)$.

     \subsection{The ${\frak sl}_2$ case.} Let us investigate
throughly quasi-Verma modules in the case of
$\U_\ell=\U_\ell({\frak sl}_2)$. First we find the simple
subqoutient modules in the module
$M_\ell^e(k\ell)=M_\ell(k\ell).$

     Recall the classification of the simple objects in the
category of $X$-graded $\U_\ell^0$-semisimple $\U_\ell$-modules
locally finite with respect to the action of $E_i$ and
$E_i^{(\ell)}$, $i\in I$, obtained by Lusztig in \cite{L2}. In
the ${\frak sl}_2$ case it looks as follows. Identify the
weight lattice $X$ with $\Z$. \vskip 1mm \noindent \Prop
\begin{itemize} \item[(i)] For $0\le k<\ell$ the simple
$\u_\ell({\frak sl}_2)$-module $L(k)$ is a restriction of a
simple $\U_\ell({\frak sl}_2)$-module. \item[(ii)] Any simple
$\U_\ell({\frak sl}_2)$-module from the category described
above is isomorphic to a module of the form $L(k)\ten
L(m\ell)$, where $0\le k<\ell$. Here the simple module
$L(m\ell)$ is obtained from the simple $U({\frak sl}_2)$-module
$L(m)$ via restriciton using the map $\U_\ell({\frak
sl}_2)\map\U({\frak sl}_2)//\u_\ell({\frak sl}_2)=U({\frak
sl}_2)$. \qed \end{itemize} Denote the only reflection in the
Weyl group for ${\frak sl}_2$ by $s$. \vskip 1mm \noindent
\Lemma \begin{itemize} \item[(i)] $\ch M_\ell^e(0)=\ch L(0)+\ch
L(-2)+\ch L(-2\ell)$. \item[(ii)] For $k>0$ we have $\ch
M_\ell^e(k\ell)=\ch L(k\ell)+\ch L(k\ell-2)+\ch L(-k\ell-2)+\ch
L(-(k+2)\ell)$. \item[(iii)] $\ch M_\ell^s(s\cdot 0)=\ch
L(-2)+\ch L(-2\ell)$. \item[(iv)] For $k>0$ we have $\ch
M_\ell^s(s\cdot k\ell)=+\ch L(-k\ell-2)+\ch L(-(k+2)\ell)$.
\qed \end{itemize} In fact it is easy to find the filtrations
on quasi-Verma modules with simple subquotients that correspond
to the character equalities above. \vskip 1mm \noindent \Lemma
\begin{itemize} \item[(i)] For $k>0$ there exist exact
sequences \begin{gather*} 0\map L(k\ell-2)\map W(k\ell)\map
L(k\ell)\map 0,\\ 0\map L(-k\ell-2)\map M_\ell^s(s\cdot
k\ell)\map L(-(k+2)\ell)\map 0. \end{gather*} \item[(ii)] There
exists a filtration on $M_\ell^e(0)$ with subquotients as
follows: $$ \gr^1M_\ell^e(0)=L(0),\ \gr^2M_\ell^e(0)=L(-2),\
\gr^3M_\ell^e(0)=L(-2\ell). $$ \item[(iii)] For $k>0$ there
exists a filtration on $M_\ell^e(k\ell)$ with subquotients as
follows: \begin{gather*} \gr^1M_\ell^e(k\ell)=L(k\ell),\
\gr^2M_\ell^e(k\ell)=L(k\ell-2),\\
\gr^3M_\ell^e(k\ell)=L(-(k+2)\ell),\
\gr^4M_\ell^e(k\ell)=L(-k\ell-2).\qed \end{gather*}
\end{itemize} Thus we obtain the following statement.

     \sssn \label{sl2case}
     \Prop For any $k\ge0$ there exists an exact sequence of
$\U_\ell=\U_\ell({\frak sl}_2)$-modules $$ 0\map
M_\ell^s(s\cdot k\ell)\map M_\ell^e(k\ell)\map
W(k\ell)\map0.\qed $$ \Rem \begin{itemize} \item[(i)]
     In fact it is easy to verify that $M_\ell^s(s\cdot k\ell)$
is isomorphic to the {\em contragradient} Verma module $\D
M_\ell(s\cdot k\ell)$.
     \item[(ii)] Note that if $\xi$ is not a root of unity then
the usual BGG resolution in the ${\frak sl}_2$ case provides an
exact complex $$ 0\map M_\xi^s(s\cdot k\ell)\map
M_\xi^e(k\ell)\map L(k\ell)\map0. $$ \end{itemize} Thus we see
that the flat family of such complexes over
$\CC^*\setminus\{$roots of unity$\}$ is extended over the whole
$\CC^*$.

     \subsection{Construction of the quasi-BGG complex.} Here
we extend the previous considerations to the case of the
quantum group $\U_\ell$ for arbitrary root data $(Y,X,\ldots)$
of the finite type $(I,\cdot)$. Fix a dominant weight $\lambda
\in X$.

     First we construct an inclusion
$M_\ell^{w'}(w'\cdot\ell\lambda)\hookrightarrow
M_\ell^{w}(w\cdot\ell\lambda)$ for a pair of elements $w',w\in
W$ such that $l(w')=l(w)+1$ and $w'>w$ in the Bruhat order on
the Weyl group. In fact we can do it explicitly only for $w'$
and $w$ differing by a simple reflection: $w'=ws_i$, $i\in I$.

     Consider the twisted quantum parabolic subalgebra
$w(\P_{i,\ell})$. Then $w(\P_{i,\ell})\supset w(\B_\ell^+)$ and
$w(\P_{i,\ell})\supset ws_i(\B_\ell^+)$. Consider also the Levi
quotient algebra $w(\P_{i,\ell})\map w(\LL_{i,\ell})$. The
algebra $w(\LL_{i,\ell})$ is isomorphic to $\U_\ell({\frak
sl}_2)\ten_{\U_\ell^0({\frak sl}_2)}\U_\ell^0$.

     By Proposition \ref{sl2case} we have a natural inclusion
of $w (\LL_{i,\ell})$-modules
$\SInd_{ws_i(\LL_{i,\ell}^+)}^{w(\LL_{i,\ell})}\CC(\ell\lambda)
\hookrightarrow
\SInd_{w(\LL_{i,\ell}^+)}^{w(\LL_{i,\ell})}\CC(\ell\lambda)$.
\vskip 1mm \noindent \Lemma \begin{itemize} \item[(i)]
$\SInd_{w(\B_\ell^+)}^{\U_\ell}(\cdot)=
\SInd_{w(\P_{i,\ell})}^{\U_\ell}\circ \SInd_{w(\B_\ell^+)}^{
w(\P_{i,\ell})}(\cdot)$. \item[(ii)]
$\SInd_{ws_i(\B_\ell^+)}^{\U_\ell}(\cdot)=
\SInd_{w(\P_{i,\ell})}^{\U_\ell}\circ
\SInd_{ws_i(\B_\ell^+)}^{w(\P_{i,\ell})}(\cdot)$. \item[(iii)]
$\SInd_{w(\B_\ell^+)}^{\U_\ell}(\CC(\ell\lambda))=
\SInd_{w(\P_{i,\ell})}^{\U_\ell}\circ
\Res_{w(\P_{i,\ell})}^{w(\LL_{i,\ell})}\circ
\SInd_{W(\LL_{i,\ell}^+)}^{w(\LL_{i,\ell})}(\ell\lambda)$.
\item[(iv)]
$\SInd_{ws_i(\B_\ell^+)}^{\U_\ell}(\CC(\ell\lambda))=
\SInd_{w(\P_{i,\ell})}^{\U_\ell}\circ
\Res_{w(\P_{i,\ell})}^{w(\LL_{i,\ell})}\circ
\SInd_{ws_i(\LL_{i,\ell}^+)}^{w(\LL_{i,\ell})}(\ell\lambda)$.
\qed \end{itemize} \noindent \Cor For $w'=ws_i>w$ in the Bruhat
order we have a natural inclusion of $\U_\ell$-modules
$i_\ell^{ws_i,w}:\
M_\ell^{ws_i}(ws_i\cdot\ell\lambda)\hookrightarrow
M_\ell^w(w\cdot\ell\lambda). \qed $

     Note that in the previous considerations we never used the
fact that $\U_\ell$ was a quantum group {\em at the root of
unity}. Thus a similar construction provides inclusions
$$i_\xi^{ws_i,w}:\
M_\ell^{ws_i}(ws_i\cdot\ell\lambda)\hookrightarrow
M_\xi^w(w\cdot\ell\lambda) $$ for any $\xi\in \CC^*$. Recall
that if $\xi$ is noot a root of unity then
$M_\xi^w(w\cdot\ell\lambda)$ is isomorphic to the usual Verma
module $M_\xi(w\cdot\ell\lambda)$. Thus the morphism
$i_\xi^{w',w}$ coincides with the standard inclusion of Verma
modules constructed by J. Bernshtein, I.M. Gelfand and S.I.
Gelfand in \cite{BGG} that becomes a component of the
differential in the BGG resolution. In other words we see that
the flat family of inclusions $i_\xi^{ws_i,w}:\
M_\ell^{ws_i}(ws_i\cdot\ell\lambda)\hookrightarrow
M_\xi^w(w\cdot\ell\lambda) $ defined over
$\CC^*\setminus\{$roots of unity$\}$ can be extended naturally
over the whole $\CC^*$.

     Iterating the inclusion maps we obtain a flat family of
submodules $i^w_\xi(M_\xi^w(w\cdot\ell\lambda))\subset
M_\xi^e(\ell\lambda)$ for $\xi\in\CC^*$, $w\in W$, providing an
extension of the standard lattice of Verma submodules in
$M_\xi(\ell\lambda)$ for $\xi\in\CC^*\setminus\{$roots of
unity$\}$.

     \sssn \Lemma For a pair of elements $w',w\in W$ such that
$l(w')=l(w)+1$ and $w'>w$ in the Bruhat order we have $$
i_\ell^{w'}(M_\ell^{w'}(w'\cdot\ell\lambda)) \hookrightarrow
i_\ell^w(M_\ell^{w}(w\cdot\ell\lambda)). $$ \dok To prove the
statement note that the condition $\{A_\xi$ is a submodule in
$B_\xi\}$ is a closed condition in a flat family. \qed

     Now using the standard combinatorics of the classical BGG
resolution we obtain the following statement.

     \sssn \Theorem There exists a complex of $\U_\ell$-modules
$B\bul_\ell(\ell\lambda)$ with $$
B_\ell^{-k}(\ell\lambda)=\underset{w\in W,l(w)=k}{\bigoplus}
M_\ell^w(w\cdot\ell\lambda) $$ and with differentials provided
by direct sums of the inclusions $i_\ell^{w',w}$. \qed

     \sssn \Def We call the complex $B\bul_\ell(\ell\lambda)$
the {\em quasi-BGG complex} for the dominant weight
$\ell\lambda\in X$.

     Recall that we have a natural morphism $can:
M_\ell^e(\ell\lambda)\map W(\ell\lambda)$. \vskip 1mm \noindent
\Prop $H^0(B\bul_\ell(\ell\lambda))=W(\ell\lambda). $

     \dok The statement follows from the standard fact that
$W(\ell\lambda)$ is the biggest finite dimensional quotient
module in $M_\ell^e(\ell\lambda)$. \qed

     \sssn \label{conexact} \Con The complex
$B\bul_\ell(\ell\lambda)$ is quasiisomorphic to the Weyl module
$W(\ell\lambda)$.\qed

     Note that by the previous subsection the Conjecture holds
in the $\U_\ell({\frak sl}_2)$ case.

     Using the contragradient duality we obtain a complex $\D
B\bul_\ell(\ell\lambda)$ consisting of direct sums of
contragradient quasi-Verma modules with $H^0(\D
B\bul_\ell(\ell\lambda))=\D W(\ell\lambda)$. In particular we
have a morphism in the category of complexes of
$\U_\ell$-modules $\D(can):\ \D W(\ell\lambda)\map\D
B_\ell\bul(\ell\lambda)$.

     \subsection{Semiinfinite cohomology with coefficients in
quasi-Verma modules.} Recall the following construction that
plays crucial role in considerations of Ginzburg and Kumar in
\cite{GK}.

     Let $\left(\B_\ell^+\mod\right)^{\operatorname{fin}}$
(resp. $\left(\U_\ell\mod\right)^{\operatorname{fin}}$, resp.
$\left(U(\b^+)\mod\right)^{\operatorname{fin}}$, resp.
$\left(U(\g)\mod\right)^{\operatorname{fin}}$) be the category
of finite dimensional $X$-graded modules over the corresponding
algebra with the action of the Cartan subalgebra semisimple and
well defined with respect to the $X$-gradings. Consider the
functors:
\begin{gather*}
\Coind_{\B_\ell^+}^{\U_\ell}:\
\B_\ell^+\mod\map\U_\ell\mod;\
\left(\Coind_{\B_\ell^+}^{\U_\ell}\right)^{\operatorname{fin}}:
\ \left(\B_\ell^+\mod\right)^{\operatorname{fin}}
\map\left(\U_\ell\mod\right)^{\operatorname{fin}};\\
     \Coind_{U(\b^+)}^{U(\g)}:\ U(\b^+)\mod\map U(\g)\mod;\
\left(\Coind_{U(\b^+)}^{U(\g)}\right)^{\operatorname{fin}}:\
\left(U(\b^+)\mod\right)^{\operatorname{fin}}\map\left(U(\g)\mod
\right)^{\operatorname{fin}},\\
     (\cdot)^{\b_\ell^+}:\ \B_\ell^+\mod\map U(\b^+)\mod \text{
and } \left(\B_\ell^+\mod\right)^{\operatorname{fin}}\map
\left(U(\b^+)\mod\right)^{\operatorname{fin}};\\
(\cdot)^{\u_\ell}:\ \U_\ell\mod\map U(\g)\mod \text{ and }
\left(\U_\ell\mod\right)^{\operatorname{fin}}\map
\left(U(\g)\mod\right)^{\operatorname{fin}},
\end{gather*}
where $(\cdot)^{\operatorname{fin}}$ denotes taking the maximal
finite dimensional submodule in $(\cdot)$ and
$(\cdot)^{\b_\ell^+}$ (resp. $(\cdot)^{\u_\ell}$) denotes
taking $\b_\ell^+$- (resp. $\u_\ell$)-invariants. \vskip 1mm
\noindent \Prop (see \cite{GK}) \begin{itemize} \item[(i)]
$(\cdot)^{\u_\ell}\circ \Coind_{\B_\ell^+}^{\U_\ell}=
\Coind_{U(\b^+)}^{U(\g)}= (\cdot)^{\b_\ell^+}$; \item[(ii)]
$(\cdot)^{\u_\ell}\circ
\left(\Coind_{\B_\ell^+}^{\U_\ell}\right)^{\operatorname{fin}}=
\left(\Coind_{U(\b^+)}^{U(\g)}\right)^{\operatorname{fin}}
(\cdot)^{\b_\ell^+}$.\qed \end{itemize} The semiinfinite
analogue for the first part of the previous statement looks as
follows. Fix $w\in W$. Consider the functors: \begin{gather*}
\SCoind_{w(\B_\ell^+)}^{\U_\ell}:\ {\mathsf
D}(w(\B_\ell^+)\mod)\map{\mathsf D} (\U_\ell\mod),\\
\SCoind_{U(w(\b^+))}^{U(\g)}:\ {\mathsf
D}(U(w(\b^+))\mod)\map{\mathsf D} (U(\g)\mod),\\ \Exts_{
w(\b_\ell^+)}(\underC,\cdot):\ {\mathsf
D}(w(\B_\ell^+)\mod)\map{\mathsf D} (U(w(\b^+))\mod),\\
\Exts_{\u_\ell}(\underC,\cdot):\ {\mathsf
D}(w(\U_\ell)\mod)\map{\mathsf D} (U(\g))\mod). \end{gather*}
\sssn \Theorem We have $$ \Exts_{\u_\ell}(\underC,\cdot) \circ
\SCoind_{w(\B_\ell^+)}^{\U_\ell} = \SCoind_{U(w(\b^+))}^{U(\g)}
\circ \Exts_{w(\b_\ell^+)}(\underC,\cdot). $$ \dok To simplify
the notations we work with semiinfinite homology and
semiinfinite induction instead of semiinfinite cohomology and
semiinfinite coinduction. By \cite{Ar1}, Appendix B, every
convex complex of $w(\B_\ell^+)$-modules is quasiisomorphic to
a K-semijective convex complex. Consider a K-semijective convex
complex of $w(\B_\ell^+)$-modules $SS\bul$. Note that both
semiinfinite induction functors are exact and take
K-semijective complexes to K-semijective complexes. Thus we
have \begin{gather*} \Tors^{\u_\ell}(\underC,\cdot)\circ
\SInd_{w(\B_\ell^+)}^{\U_\ell}(SS\bul)=
\Tors^{\u_\ell}(\underC,
\Tors^{w(\B_\ell^+)}(S_{\U_\ell}^{\U_\ell^+},SS\bul))\\=
\Tors^{w(\B_\ell^+)}( \Tors^{\u_\ell}(\underC,
S_{\U_\ell}^{\U_\ell^+}),SS\bul))= \Tors^{w(\B_\ell^+)}
(S_{U(\g)}^{U(\n^+)},SS\bul)\\= \Tors^{w(\B_\ell^+)} (\underC,
S_{U(\g)}^{U(\n^+)})\ten SS\bul)= \Tors^{U(w(\b^+))}( \underC,
\Tors^{w(\b_\ell^+)}(\underC, S_{U(\g)}^{U(\n^+)}) \ten
SS\bul))\\= \Tors^{U(w(\b^+))}( S_{U(\g)}^{U(\n^+)}),
\Tors^{w(\b_\ell^+)}( \underC,SS\bul)) =
\SInd_{U(w(\b^+))}^{U(\g)} \circ
\Tors^{w(\b_\ell^+)}(\underC,SS\bul). \end{gather*} Here we
used the fact that the subalgebra $w(\b_\ell^+)\subset
w(\B_\ell^+)$ is normal with the quotient algebra equal to
$U(w(\b^+))$. \qed \vskip 1mm \noindent \Cor $$
\Exts_{\u_\ell}(\underC,\D M_\ell^w(w\cdot\ell\lambda))=
H_{T_{S_w}^*(G/B)}^{\sharp(R^+)}(T^*(G/B),\pi^*\L(\lambda))$$
     as a module over both $U(\g)$ and
$H^0(T^*(G/B),\O_{T^*(G/B)})=H^0(\CN,\O_{\CN})$.
     \qed

     Recall that for the Springer-Grothendieck resolution of
the nilpotent cone $\mu:\ T^*(G/B)\map \CN$ we have
$\mu^{-1}(\n^+)=\underset{w\in W}{\bigsqcup}T^*_{S_w}( G/B)$.
\vskip 1mm \noindent \Prop \begin{itemize} \item[(i)] There
exists a filtration on
$H^{\sharp(R^+)}_{\mu^{-1}(\n^+)}(T^*(G/B),\pi^*\L(\lambda))$
with the subquotients equal to
$H^{\sharp(R^+)}_{T^*_{S_w}(G/B)}(T^*(G/B),\pi^*\L(\lambda))$,
for $w\in W$. \item[(ii)] $ \Exts_{\u_\ell}(\underC,\D
B_\ell\bul(\ell\lambda))=
H_{\mu^{-1}(\n^+)}^{\sharp(R^+)}(T^*(G/B),\pi^*\L(\lambda))$.
\item[(iii)] The morphism \begin{gather*}
H_{\n^+}^{\sharp(R^+)}(\CN,\mu_*\pi^*\L(\lambda))
\overset{\sigma}{\map} \Exts_{\u_\ell}(\underC,\D
W(\ell\lambda)) \overset{\tau}{\map} \Exts_{\u_\ell}(\underC,\D
B_\ell\bul(\ell\lambda))\\ \til{\map}
H_{\mu^{-1}(\n^+)}^{\sharp(R^+)}(T^*(G/B),\pi^*\L(\lambda))
\end{gather*} is an isomorphism. \qed \end{itemize} \sssn \Cor
\label{half} \begin{itemize} \item[(i)] The map $\sigma$ is
injective and the map $\tau$ is surjective. \item[(ii)]
Conjecture \ref{maincon} is true if conjecture~\ref{conexact}
is so. \qed \end{itemize} However we will manage to prove
Conjecture~\ref{maincon} without using the exactness of the
quasi-BGG complex.

     \subsection{Quantum twisting functor.} Here we present a
construction of a kernel functor similar to the trwisting
functors in \cite{Ar5} and \cite{Ar6}. Consider the subalgebra
$\U_\ell^{1/2}$ in $\U_\ell$ generated by the elements $E_i$,
$F_i$, $F_i^{(\ell)}$, $i\in I$. \vskip 1mm \noindent \Lemma
\begin{itemize} \item[(i)] The subalgebra in $\U_\ell$
generated by the elements $F_i^{(\ell)}$, $i\in I,$ is
isomorphic to $U(\n^-)$. \item[(ii)] The algebra
$\U_\ell^{1/2}$ is isomorphic to $\u_\ell\ten U(\n^-)$ as a
vector space. \item[(iii)] The subalgebra
$\u_\ell\subset\U_\ell^{1/2}$is normal with the quotient
algebra equal to $U(\n^-)$. \qed \end{itemize} Note yet that
this splitting of the projection $\U_\ell^{1/2}\map U(\n^-)$
does not extend to a splitting of the projection $\U_\ell\map
U(\g)$.

     Consider the restriction of the quasi-Verma module
$M_\ell^w(w\cdot\ell\lambda)$ onto the subalgebra
$\U_\ell^{1/2}$. \vskip 1mm

     \noindent \Lemma The $\U_\ell^{1/2}$-module
$M_\ell^w(w\cdot\ell\lambda)$ is isomorphic to $
\Ind_{w(\b_\ell^+)}^{\U_\ell^{1/2}}\CC(\ell(w\cdot\lambda))
$.\qed

     In particular $M_\ell^w(w\cdot\ell\lambda)$ is {\em free}
over the subalgebra $U(\n^-)\subset \U_\ell$. Consider the left
semiregular $\U_\ell$-module $S_{\U_\ell}^{\U_\ell^-}$. \vskip
1mm \noindent \Lemma
$S_{\U_\ell}^{\U_\ell^-}\til{\map}\Ind_{U(\n^-)}^{\U_\ell}U(\n^
-)^*$ as a left
     $\U_\ell$-module. \qed \vskip 1mm \noindent \Cor
\begin{itemize} \item[(i)]
$\left(S_{\U_\ell}^{\U_\ell^-}\ten_{\U_\ell}M_\ell^w(w\cdot\ell
\lambda)\right)^\tau
     = \D M_\ell^{w_0w^{-1}}(w_0w^{-1}\cdot\ell\lambda)$.
\item[(ii)]
$\left(S_{\U_\ell}^{\U_\ell^-}\overset{\operatorname{L}}{\ten}_
{\U_\ell}W(\ell\lambda)\right)^\tau
     = \D W(\ell\lambda)[\sharp(R^+)]$. \item[(iii)]
$\left(S_{\U_\ell}^{\U_\ell^-}\ten_{\U_\ell}B_\ell\bul(\ell\lambda)
\right)^\tau
     = \D B_\ell\bul(\ell\lambda)[\sharp(R^+)]$. Here
$(\cdot)^\tau$ denotes twisting the $\U_\ell$-action on
$(\cdot)$ by the Chevalley involution $\tau$. \qed
\end{itemize} We denote the functor
$\left(S_{\U_\ell}^{\U_\ell^-}\overset{\operatorname{L}}{\ten}_
{\U_\ell}\cdot\right)^\tau$
     by ${\frak S}_\ell$. Now consider the morphism of
complexes of $\U_\ell$-modules $can:\
B_\ell\bul(\ell\lambda)\map W(\ell\lambda)$. We have a morphism
in the derived category of $\U_\ell$-modules ${\frak
S}_\ell(can):\ \D B_\ell\bul(\lambda)[\sharp(R^+)]\map \D
W(\ell\lambda)[\sharp(R^+)]$. On the other hand consider the
canonical map $\D(can):\ \D W(\ell\lambda) \map\D
B_\ell\bul(\lambda)$.

     \sssn \Prop For $\ell$ big enough we have ${\frak
S}_\ell(can)\circ \D(can)=\operatorname{Id}_{\D
W(\ell\lambda)[\sharp(R^+)]}$.

     \dok Note that the in the definition of the functor
${\frak S}_\ell$ we never used the fact that we worked in the
derived category of modules over the quantum group {\em at a
root of unity}. In particular one can construct similar
functors $$ {\frak S}_\xi:\ {\mathsf D}(\U_\xi\mod)\map{\mathsf
D}(\U_\xi\mod) $$ for any $\xi\in \CC^*$. It is easy to verify
that for $\xi$ not a root of unity we have ${\frak
S}_\xi(can)\circ
\D(can)=\operatorname{Id}_{L(\ell\lambda)[\sharp(R^+)]}$. On
the other hand note that the endomorphism spaces of the
$\U_\xi$-modules $L(\ell\lambda)$ (resp. $\U_\ell$-modules $\D
W(\ell\lambda)$) is one dimensional. Thus we can consider
${\frak S}_\xi(can)\circ \D(can)$ as a (nonzero) polynomial
function on $\xi$. But the number of roots of this polynomial
is finite. \qed \vskip 1mm \noindent \Cor
$\Exts_{\u_\ell}(\underC,\D W(\ell\lambda)) $ is a direct
summand in $\Exts_{\u_\ell}(\underC,\D B_\ell\bul(\ell\lambda))
$. \qed

     Comparing this statement with Corollary~\ref{half} we
obtain the main result of the section.

     \sssn \Theorem For $\ell$ big enough
$\ext_{\u_\ell}^{\si+\bullet}(\underC,\D W(\ell\lambda))=
H^0_{\n^+}(\CN,\mu_*p^*\L(\lambda))$
     as a $H^0(\CN,\O_{\CN})$-module. \qed

     \section{Further results and conjectures.} In this section
we present several facts without proof. We also formulate some
conjectures concerning possible origin of the quasi-BGG
complex.

     \subsection{Alternative triangular decompositions of
$\u_\ell$.} \label{alt} Note that the definition of
semiinfinite cohomology starts with specifying a triangular
decomposition of a {\em graded} algebra $\u$. Fix a subset
$J\subset I$. Instead of the usual height function consider the
linear map $\hgt_J:\ X\map\Z$ defined on the elements $i', i\in
I$ by $\hgt_J(i')=0$ for $i\in J$ and $\hgt_J(i')=1$ otherwize
and extended to the whole $X$ by linearity. Now we work in the
category of complexes of $X$-graded $\u_\ell$-modules
satisfying conditions of concavity and convexity with respect
to the $\Z$-grading obtained from the $X$-grading with the help
of the function $\hgt_J$.

     Consider the triangular decomposition of the small quantum
group $\u_\ell=\p_{J,\ell}^-\ten \u_{J,\ell}^+$, where
$\p_{J,\ell}^-$ denotes the small quantum negative parabolic
subalgebra in $\u_\ell$ corresponding to the subset $J\subset
I$ and $\u_{J,\ell}^+$ denotes the quantum analogue

     of the nilpotent radical in $\p_J^+$ defined with the help
of Lusztig generators of $\U_\ell$ and $\u_\ell$ (see
\cite{L1}).

     Then it is known that the subalgebra $\u_{J,\ell}^+$ in
$\u_\ell$ is Frobenius just like $\u_\ell^+$. Thus it is
possible to use the general definition of semiinfinite
cohomology presented in \ref{setup}. Denote the corresponding
functor by $\Exts_{\u_\ell,J}(*,*)$.

     On the other hand consider the classical negative
parabolic subalgebra $\p_J^-\subset\g$ and its nillradical
$\n_J^-$. Choose the standard $X$-homogeneous root basis
$\{f_\alpha\}$ in the space $\n_J^-$. Consider the subset in
$\CN^{(J)}\subset\CN$ annihilated by all the elements of the
base dual to $\{f_\alpha\}$.

     \sssn \Theorem
$\Exts_{\u_\ell,J}(\underC,\underC)\til{\map}
H^{\dim(\n^-_J)}_{\CN^{(J)}}(\CN,\O_{\CN})$ as
$H^0(\CN,\O_{\CN})$-modules. \qed

     \subsection{Contragradient Weyl modules with
non-$\ell$-divisible highest weights.} Fix a {\em dominant}
weight of the form $\ell\lambda+ w\cdot 0$, where $\lambda\in
X$ and $w\in W$. It is known that all the dominant weights in
the linkage class containing $0$ look like this. Consider the
contragradient Weyl module $\D W(\ell\lambda+ w\cdot0)$. The
following statement generalizes Corollary \ref{chformula}. Its
proof is similar to the proof of Conjecture \ref{maincon}.

     \sssn \Theorem \begin{gather*} \ch\left(
\ext_{\u_\ell}^{\si+\bullet}(\underC,\D W(\ell\lambda+
w\cdot0)),t\right)\\=
\frac{t^{-\sharp(R^+)+l(w)}}{\prod_{\alpha\in
R^+}(1-e^{-\ell\alpha})}\sum_{v\in W}
\frac{e^{v(\ell\lambda)}t^{2l(v)}}{\prod_{\alpha\in
R^+,v(\alpha)\in R^+}(1-t^2e^{-\ell\alpha}) \prod_{\alpha\in
R^+,v(\alpha)\in R^-}(1-t^{-2}e^{-\ell\alpha})}.\qed
\end{gather*} \subsection{Contragradient Weyl modules:
alternative triangular decompositions.} Fix the triangular
decomposition of the small quantum group $\u_\ell$ like in
\ref{alt}. A natural generalization of Conjecture~\ref{maincon}
to the case of the parabolic triangular decomposition looks as
follows. We keep the notations from~\ref{alt}.

     \sssn \Con $\Exts_{\u_\ell,J}(\underC,\D
W(\ell\lambda))\til{\map}
H^{\dim(\n^-_J)}_{\CN^{(J)}}(\CN,\mu_*p^*\L(\lambda))$ as a
$H^0(\CN,\O_{\CN})$-module. \qed

     \subsection{Connection with affine Kac-Moody algebras}
Finally we would like to say a few words about a possible
explanation for the existence of quasi-Verma modules and
quasi-BGG resolutions.

     Suppose for simplicity that the root data $(Y,X,\ldots)$
are {\em untwisted}, i. e. the corresponding Cartan matrix is
symmetric. Consider the affine Lie algebra
$\hat\g=\g\ten\CC[t,t^{-1}]\oplus\CC K$ corresponding to $\g$.
Fix a {\em negative} level $-h^\vee+k$, where $k\in
1/2\Z_{<0}\setminus\Z_{<0}$ and $h^\vee$ denotes the dual
Coxeter number for chosen root data of the finite type.
Consider the Kazhdan-Lusztig category $\til{\mathcal O}_{-k}$
of $\g\ten\CC[t]$-integrable finitely generated
$\hat\g$-modules diagonizible with respect to the Cartan
subalgebra in $\hat\g$ at the level $-2h^\vee-k$. Kazhdan and
Lusztig showed that the category $\til{\mathcal O}_k$ posasses
a structure of a rigid tensor category with the {\em fusion}
tensor product $\overset{\cdot}{\ten}$. Moreover, they proved
the following statement.

     \sssn \Theorem (see \cite{KL}) Let $\ell=-2k$. Yhen the
tensor category $(\til{\mathcal O}_k,\overset{\cdot}{\ten})$ is
equivalent to the category
$\left(\U_\ell\mod\right)^{\operatorname{fin}}$ with the tensor
product provided by the Hopf algebra structure on $\U_\ell$.
\qed

     We denote the functor $ (\til{\mathcal
O}_k,\overset{\cdot}{\ten})\map
\left(\left(\U_\ell\mod\right)^{\operatorname{fin}},\ten\right)
$
     providing the equivalence of categories by $\til{\mathsf
{kl}}$. Consider the usual category ${\mathcal O}_k$ for
$\hat\g$ at the same level. Finkelberg constructed a functor
${\mathsf {kl}}:\ {\mathcal O}_k\map\U_\ell\mod$ extending the
functor $\til{\mathsf {kl}}$ (see \cite{F}). Note that the
functor $\mathsf{kl}$ has no chance to be an equivalence of
categories because it is known not to be exact.

     Fix a dominant (resp. arbitrary) weight $\lambda\in X$.
     Consider now the contragradient Weyl module $\D{\mathcal
W}(\lambda)=\Coind_{U(\g\ten\CC[t])}^{U_{-2h^\vee-k}(\hat\g)}L(
\lambda)$
     and the contragradient Verma module $\D{\mathcal
M}(\lambda)=\Coind_{U(\g\ten\CC[t])}^{U_{-2h^\vee-k}(\hat\g)}\D
M(\lambda)$ over $\hat\g$ at the level $-2h^\vee-k$, where
$L(\lambda)$ (resp. $M(\lambda)$) denotes the simple module
(resp. the contragradient Verma module over $\g$ with the
highest weight $\lambda$. Then the usual contragradient BGG
resolution of $L(\lambda)$ provides a resolution
$\D\B(\lambda)$ of the contragradient Weyl module $\D {\mathcal
W}(\lambda)$ consisting of direst sums of contragradient Verma
modules of the form $\D{\mathcal M}(w\cdot\lambda)$, where $\in
W$. It is known that the Kazhdan-Lusztig functor takes Weyl and
contragradient Weyl modules over $\hat\g$ to Weyl (resp.
contragradient Weyl) modules over $\U_\ell$.

     \sssn \Con The functor ${\mathsf {kl}}$ takes $\D{\mathcal
M}(w\cdot\ell\lambda)$ to the contragradient quasi-Verma module
$\D M_\ell^w(w\cdot\ell\lambda)$. Moreover the complex
${\mathsf {kl}}(\D{\mathcal B}(\ell\lambda))$ is
quasiisomorphic to $\D W(\ell\lambda)$.\qed

\end{document}